\documentclass{article}

\title{\LARGE \textbf{Long cycles in graphs through fragments }}
\author{Zh.G. Nikoghosyan\footnote{G.G. Nicoghossian (up to 1997)}\\
Institute for Informatics and Automation Problems\\ National Academy of Sciences\\
P.Sevak 1, Yerevan 0014, Armenia\\ E-mail: zhora@ipia.sci.am}

\begin{document}

\maketitle

\begin{abstract}
Four basic Dirac-type sufficient conditions for a graph $G$ to be hamiltonian  are known involving order $n$, minimum degree $\delta$, connectivity $\kappa$ and independence number $\alpha$ of $G$: (1) $\delta \geq n/2$ (Dirac); (2) $\kappa \geq 2$ and $\delta \geq (n+\kappa)/3$ (by the author); (3) $\kappa \geq 2$ and $\delta \geq \max\lbrace (n+2)/3,\alpha \rbrace$ (Nash-Williams); (4) $\kappa \geq 3$ and $\delta \geq \max\lbrace (n+2\kappa)/4,\alpha \rbrace$ (by the author). In this paper we prove the reverse version of (4) concerning the circumference $c$ of $G$ and completing the list of reverse versions of (1)-(4): (R1) if $\kappa \geq 2$, then $c\geq\min\lbrace n,2\delta\rbrace$ (Dirac); (R2) if $\kappa \geq 3$, then $c\geq\min\lbrace n,3\delta -\kappa\rbrace$ (by the author); (R3) if $\kappa\geq 3$ and $\delta\geq \alpha$, then $c\geq\min\lbrace n,3\delta-3\rbrace$ (Voss and Zuluaga); (R4) if $\kappa\geq 4$ and $\delta\geq \alpha$, then $c\geq\min\lbrace n,4\delta-2\kappa\rbrace$. To prove (R4), we present four more general results centered around a lower bound $c\geq 4\delta-2\kappa$ under four alternative conditions in terms of fragments. A subset $X$ of $V(G)$ is called a fragment of $G$ if $N(X)$ is a minimum cut-set and $V(G)-(X\cup N(X))\neq\emptyset$. \\
  

\end{abstract}

\section{Introduction}

The classic hamiltonian problem asks to check whether a given graph has a spanning cycle. Such cycles are called Hamilton cycles in honor of Sir William Rowan Hamilton, who, in 1856, described an idea for a game. The hamiltonian problem is based entirely on two genuine concepts "graph" and "Hamilton cycle". Since this problem is NP-complete, generally it is senseless to expect nontrivial results in this area within these two initial concepts and it is natural to look for conditions for the existence of a Hamilton cycle either involving quite new concepts or transforming the initial ones. 

In 1952, Dirac [2] obtained the first sufficient condition for a graph to be hamiltonian based on "minimum degree $\delta$". Actually, this successful combination of three genuine concepts "graph", "Hamilton cycle" and "minimum degree" marked the beginning of a new period in hamiltonism generating a wide class of various problems and ideas for fruitful explorations. Further, these concepts continually were transformed in a way of various limitations, generalizations, extensions and manipulations based on: 

(i) structural limitations on graphs: regular and bipartite graphs, graphs with forbidden subgraphs (for example, claw-free graphs and planar graphs) and so on,

(ii) quantitative limitations (relations) on graphs: 2-connected graphs, 1-tough graphs, graphs with $\delta\geq n/2$ and so on,

(iii) generalized Hamilton cycles: long cycles, Hamilton paths and their generalizations (for example, spanning trees with minimum number of leaves), 2-factors, large cycles (for example, dominating cycles and generalized dominating cycles with complements of certain structures) and so on,

(iv) generalized minimum degree notions: degree sequences, degree sums, neighborhood unions, generalized degrees and so on.

Due to transformations (i)-(iv), the frames of a concept "hamiltonian problem" were expanded rapidly involving various related concepts and occupying the major directions in so called "hamiltonian graph theory".

As for minimum degree (Dirac-type) approach, it has been inspired by a couple of well-known results (direct and reverse versions) due to Dirac [2], determining how small the minimum degree $\delta$ of a graph $G$ must be to guarantee the existence of a Hamilton cycle and how large is the circumference $c$ (the length of a longest cycle) depending on $\delta$. Although the corresponding starting bounds $n/2$ and $\min \{n,2\delta\}$ in these theorems are best possible, since 1952 a number of other analogous best possible theorems appeared essentially lowering the bound $n/2$ and enlarging the bound $2\delta$ due to direct incorporation of some additional graph invariants into these bounds. 

At present, four basic Dirac-type hamiltonian sufficient conditions are known directly involving order $n$, minimum degree $\delta$, connectivity $\kappa$ and independence number $\alpha$ with minimum additional limitations and transformations of the initial conceptions due to Dirac [2], the author [8],[9], Nash-Williams [7] and the author [10], respectively.\\

\noindent\textbf{Theorem A} [2]. Every graph with $\delta\geq\frac{1}{2}n$ is hamiltonian.\\

\noindent\textbf{Theorem B} [9]. Every 2-connected graph with $\delta \geq \frac{1}{3}(n+\kappa)$ is hamiltonian.\\

\noindent\textbf{Theorem C} [7]. Every 2-connected graph with $\delta \geq \max\{\frac{1}{3}(n+2),\alpha\}$ is hamiltonian.\\

\noindent\textbf{Theorem D} [10]. Every 3-connected graph with $\delta\geq\max\lbrace \frac{1}{4}(n+2\kappa),\alpha\rbrace$ is hamiltonian.\\

A short proof of Theorem B was given in [3] due to H\"{a}ggkvist. 

The reverse versions of Theorems A-C concerning long cycles in graphs, are due to Dirac [2], the author [8],[9] and Voss and Zuluaga [14], respectively. In this paper we present the detailed proof of the last reverse version corresponding to Theorem D (it was announced still in 1985 with a short outline of the proof [11]) completing the list of reverse versions of Theorems A-D.\\

\noindent\textbf{Theorem E} [2]. Every 2-connected graph has a cycle of length at least $\min\lbrace n,2\delta\rbrace$.\\

\noindent\textbf{Theorem F} [9]. Every 3-connected graph has a cycle of length at least $\min\lbrace n,3\delta -\kappa\rbrace$.\\

\noindent\textbf{Theorem G} [14]. Every 3-connected graph with $\delta \geq \alpha$ has a cycle of length at least $\min \lbrace n, 3\delta -3\rbrace$.\\

\noindent\textbf{Theorem 1} [11]. Every 4-connected graph with $\delta\geq\alpha$ has a cycle of length at least $\min\lbrace n,4\delta-2\kappa\rbrace$.\\

To prove Theorem 1, we present four more general Dirac-type results centered around a lower bound $c\geq4\delta-2\kappa$ under four alternative conditions in terms of fragments. 

If $X\subset V(G)$, then $N(X)$ denotes the set of all vertices of $G-X$ adjacent to vertices in $X$. Furthermore, $\hat{X}$ is defined as $V(G)-(X\cup N(X))$. Following Hamidoune [6], we define a subset $X$ of $V(G)$ to be a fragment of $G$ if $N(X)$ is a minimum cut-set and $\hat{X}\neq\emptyset$. If $X$ is a fragment then $\hat{X}$ is a fragment too and $\hat{\hat{X}}=X$. For convenience, we will use $X^{\uparrow}$ and $X^{\downarrow}$ to denote $X$ and $\hat{X}$, respectively. An endfragment is a fragment that contains no other fragments as a proper subset.\\

\noindent\textbf{Theorem 2}. Let $G$ be a 3-connected graph with $\delta\geq\alpha$. If $|A^{\uparrow}|\leq 3\delta-\kappa-4$ and $|A^{\downarrow}|\leq 3\delta-3\kappa$ for an endfragment $A^{\downarrow}$ of $G$, then $c\geq\min\lbrace n,4\delta-2\kappa\rbrace$.\\

\noindent\textbf{Theorem 3}. Let $G$ be a 4-connected graph with $\delta\geq\alpha$. If $|A^{\uparrow}|\leq 3\delta-\kappa-4$, $|A^{\downarrow}|\geq 3\delta-3\kappa+1$ and $|A^{\uparrow}|\geq|A^{\downarrow}|$ for an endfragment $A^{\downarrow}$ of $G$, then $c\geq\min\lbrace n,4\delta-2\kappa\rbrace$.\\

\noindent\textbf{Theorem 4}. Let $G$ be a 4-connected graph with $\delta\geq\alpha$. If $|A^{\uparrow}|\geq 3\delta-\kappa-3$ and $|A^{\downarrow}|\leq 3\delta-3\kappa$ for an endfragment $A^{\downarrow}$ of $G$, then $c\geq\min\lbrace n,4\delta-2\kappa\rbrace$.\\

\noindent\textbf{Theorem 5}. Let $G$ be a 4-connected graph with $\delta\geq\alpha$. If $|A^{\uparrow}|\geq 3\delta-\kappa-3$ and $|A^{\downarrow}|\geq 3\delta-3\kappa+1$ for an endfragment $A^{\downarrow}$ of $G$, then $c\geq\min\lbrace n,4\delta-2\kappa\rbrace$.\\

Observe that the bounds $n/2$ and $\min\{n,2\delta\}$ in Theorems A and E were improved to $(n+\kappa)/3$ and $\min\{n,3\delta-\kappa\}$ (Theorems B and F), respectively, by direct incorporation of connectivity $\kappa$ into these bounds. We conjecture that the last two bounds are the best in a sense that they cannot be improved by an analogous way within graph invariants determinable in polynomial time. \\

\noindent\textbf{Conjecture 1}. The bounds $\frac{1}{3}(n+\kappa)$ and $3\delta-\kappa$ in Theorems B and F, respectively, can not be improved by direct incorporation of any graph invariants determinable in polynomial time.\\

\section{Definitions and notations}

By a graph we always mean a finite undirected graph $G$ without loops or multiple edges. A good reference for any undefined terms is [1]. For $H$ a subgraph of $G$ we will denote the vertices of $H$ by $V(H)$ and the edges of $H$ by $E(H)$. For every $S\subset V(G)$ we use $G-S$ short for $\langle V(G)-S\rangle$, the subgraph of $G$ induced by $V(G)-S$. In addition, for a subgraph $H$ of $G$ we use $G-H$ short for $G-V(H)$. If $X\subseteq V(G)$, then $N(X)$ denotes the set of all vertices of $G-X$ adjacent to vertices in $X$. 

Let $\delta$ denote the minimum degree of vertices of $G$. The connectivity $\kappa$ of $G$ is the minimum number of vertices whose removal from $G$ results in a disconnected or trivial graph. We say that $G$ is $s$-connected if $\kappa\ge s$. A set $S$ of vertices is independent if no two elements of $S$ are adjacent in $G$. The cardinality of maximum set of independent vertices is called the independence number and denoted by $\alpha$. 

Paths and cycles in a graph $G$ are considered as subgraphs of $G$. If $Q$ is a path or a cycle of $G$, then the length of $Q$, denoted by $|Q|$, is $|E(Q)|$. Throughout the paper the vertices and edges of a graph can be interpreted as cycles of lengths 1 and 2, respectively. A graph $G$ is hamiltonian if it contains a Hamilton cycle (a cycle containing every vertex of $G$). 

Let $C$ be a cycle of $G$ with a fixed cyclic direction. In that context, the $h$-th successor and the $h$-th predecessor of a vertex $u$ on $C$ are denoted by $u^{+h}$ and $u^{-h}$, respectively. If $h=1$, we abbreviate $u^{+1}$ and $u^{-1}$ to $u^{+}$ and $u^{-}$, respectively. For a subset $S$ of $V(C)$, we define $S^{+}=\lbrace u^{+}|u\in S\rbrace$ and $S^{-}=\lbrace u^{-}|u\in S\rbrace$. For two vertices $u$ and $v$ of $C$, let $u\overrightarrow{C}v$ denote the segment of $C$ from $u$ to $v$ in the chosen direction on $C$ and $u\overleftarrow{C}v$ denote the segment in the reverse direction. We also use similar notation for a path $P$ of $G$. For $P$ a path of $G$, denote by $F(P)$ and $L(P)$ the first and the last vertices of $P$, respectively.   

Let $Q$ be a cycle or a path of a graph $G$, $r\geq 2$ a positive integer and $Z_1,Z_2,...,Z_p$ are subsets of $V(Q)$ with $p\geq 2$. A collection $(Z_{1},...,Z_{p})$ is called a $(Q,r)$-scheme if $d_{Q}(x,y)\geq 2$ for each distinct $x,y\in Z_{i}$ (where $i\in\lbrace1,...,p\rbrace)$ and $d_{Q}(x,y)\geq r$ for each distinct $x\in Z_{i}$ and $y\in Z_{j}$ (where $i,j\in\{1,...,p\}$ and $i\neq j)$. A $(Q,r)$-scheme is nontrivial if $(Z_{1},...,Z_{p})$ has a system of distinct representatives. The definition of $(Q,r)$-scheme was first introduced by Nash-Williams [7] for $p=2$. 

Given four integers $a,b,t,\kappa$ with $\kappa\leq t$, we will use $H(a,b,t,\kappa)$ as a limit example for Theorem 1 to denote the graph obtained from $tK_a+\overline{K}_t$ by taking any $\kappa$ vertices in subgraph $\overline{K}_t$ and joining each of them to all vertices of $K_b$.\\

\noindent\textbf{Definition A $\{Q^{\uparrow}_{1},...,Q^{\uparrow}_{m};Q^{\uparrow};V^{\uparrow}_{1},...,V^{\uparrow}_{m};V^{\uparrow}\}$}. Let $A^{\uparrow}$ be a fragment of $G$ with respect to a minimum cut-set $S$. Define $Q^{\uparrow}_{1},...,Q^{\uparrow}_{m}$ as a collection of vertex disjoint paths in $\langle A^{\uparrow}\cup S\rangle$ with terminal vertices in $S$ such that $|V(Q^{\uparrow}_i)|\geq 2$ $(i=1,...,m)$ and $\sum_{i=1}^{m}|V(Q^{\uparrow}_{i})|$ is as great as possible. Abbreviate $V^{\uparrow}_{i}=V(Q^{\uparrow}_{i})$ $(i=1,...,m)$ and $V^{\uparrow}=\bigcup^{m}_{i=1}V^{\uparrow}_{i}$. Form a united path $Q^{\uparrow}$ with vertex set $V^{\uparrow}$ consisting of $Q^{\uparrow}_{1},...,Q^{\uparrow}_{m}$ and some appropriate extra-edges added in $\langle S\rangle$.\\ 

\noindent\textbf{Definition B $\{Q^{\downarrow}_{1},...,Q^{\downarrow}_{m};Q^{\downarrow}_{0};V^{\downarrow}_{1},...,V^{\downarrow}_{m};V^{\downarrow}\}$}. Let $A^{\uparrow}$ be a fragment of $G$ with respect to a minimum cut-set $S$ and $Q^{\uparrow}_{1},...,Q^{\uparrow}_{m}$ are as defined in Definition A. Denote by $Q^{\downarrow}_{1},...,Q^{\downarrow}_{m}$ a collection of paths (if exist) in $\langle A^{\downarrow}\cup S\rangle$ with $\sum_{i=1}^{m}|V(Q^{\downarrow}_{i})|$ as great as possible such that combining $Q^{\uparrow}_{1},...,Q^{\uparrow}_{m}$ with $Q^{\downarrow}_{1},...,Q^{\downarrow}_{m}$ results in a simple cycle. Abbreviate $V^{\downarrow}_{i}=V(Q^{\downarrow}_{i})$ $(i=1,...,m)$ and $V^{\downarrow}=\bigcup^{m}_{i=1}V^{\downarrow}_{i}$. For the special case $|V^{\downarrow}\cap S|=2$ and $S-V^{\uparrow}\neq\emptyset$, say $z\in S-V^{\uparrow}$, we will use $Q^{\downarrow}_{0}$ to denote a longest path in $\langle A^{\downarrow}\cup \{F(Q^{\uparrow}_{1}),L(Q^{\uparrow}_{1}),z\}\rangle$ connecting $F(Q^{\uparrow}_{1})$ and $L(Q^{\uparrow}_{1})$ and passing through z. \\ 

\noindent\textbf{Definition C $\{C^{*};C^{**}\}$}. Denote by $C^{*}$ the cycle (if exist) consisting of $Q^{\uparrow}_{1},...,Q^{\uparrow}_{m}$ and $Q^{\downarrow}_{1},...,Q^{\downarrow}_{m}$. Assume w.l.o.g. that $Q^{\uparrow}_{1},...,Q^{\uparrow}_{m}$ is chosen such that $C^{*}$ has a maximal length.  Denote by $C^{**}$ a longest cycle of $G$ with $V(C^{*})\subseteq V(C^{**})$.\\

\section{Preliminaries}

In [7], Nash-Williams proved the following result concerning $(C,r)$-schemes for a cycle $C$ and a pair $(Z_1 ,Z_2 )$ of subsets of $V(C)$.\\

\noindent\textbf{Lemma A [7]}. Let $C$ be a cycle and $(Z_{1},Z_{2})$ be a nontrivial $(C,r)$-scheme. Then

$$|V(C)|\geq\min\Big\{ 2(|Z_{1}|+|Z_{2}|)+2r-6,\frac{1}{2}r(|Z_{1}|+|Z_{2}|)\Big\}.$$

Basing on proof technique used in [7], we prove two analogous results for the families $(Z_{1},Z_{2},Z_{3})$ and $(Z_{1},Z_{2},Z_{3},Z_{4})$ of subsets of $V(C)$ under additional limitations $|Z_{1}|=1$ and $|Z_{1}|=|Z_{2}|=1$, respectively.\\

\noindent\textbf{Lemma 1}. Let $C$ be a cycle and $(Z_{1},Z_{2},Z_{3})$ be a nontrivial $(C,r)$-scheme with $|Z_{1}|=1$. Then
$$
 |V(C)|\geq\min\bigg\{ 2\sum^{3}_{i=1}|Z_{i}|+3r-12,\frac{1}{2}r\Big(\sum^{3}_{i=1}|Z_{i}|-1\Big)\bigg\}.
$$

\noindent\textbf{Lemma 2}. Let $C$ be a cycle and $(Z_{1},Z_{2},Z_{3},Z_{4})$ be a nontrivial $(C,r)$-scheme with $|Z_{1}|=|Z_{2}|=1$. Then

$$
|V(C)|\geq\min\bigg\lbrace 2\sum^{4}_{i=1}|Z_{i}|+4r-18,\frac{1}{2}r\Big(\sum^{4}_{i=1}|Z_{i}|-2\Big)\bigg\rbrace.
$$

In this paper a number of path-versions of Lemmas A, 1 and 2 will be used for the path $Q$ and the families $(Z_{1},Z_{2})$, $(Z_{1},Z_{2},Z_{3})$, $(Z_{1},Z_{2},Z_{3},Z_{4})$ of subsets of $V(Q)$ under additional limitations $|Z_{1}|=1$, $|Z_{1}|=|Z_{2}|=1$ and $|Z_{1}|=|Z_{2}|=|Z_3|=1$ in some of them.\\

\noindent\textbf{Lemma 3}. Let $Q$ be a path and $(Z_{1},Z_{2})$ be a nontrivial $(Q,r)$-scheme. Then

$$
|V(Q)|\geq\min\Big\lbrace 2(|Z_{1}|+|Z_{2}|)+r-5,\frac{1}{2}r(|Z_{1}|+|Z_{2}|-2)+1\Big\rbrace.
$$

\noindent\textbf{Lemma 4}. Let $Q$ be a path and $(Z_{1},Z_{2})$ be a nontrivial $(Q,r)$-scheme with $|Z_1|=1$ and $|Z_2|\geq2$. Then $|V(Q)|\geq2|Z_2|+r-3$.\\

\noindent\textbf{Lemma 5}. Let $Q$ be a path and $(Z_{1},Z_{2},Z_{3})$ be a nontrivial $(Q,r)$-scheme with $|Z_{1}|=1$. Then

$$
|V(Q)|\geq\min\bigg\lbrace 2\sum^{3}_{i=1}|Z_{i}|+2r-11,\frac{1}{2}r\Big(\sum^{3}_{i=1}|Z_{i}|-3\Big)+1\bigg\rbrace.
$$

\noindent\textbf{Lemma 6}. Let $Q$ be a path and $(Z_{1},Z_{2},Z_{3})$ be a nontrivial $(Q,r)$-scheme with $|Z_{1}|=|Z_2|=1$ and $|Z_3|\geq3$. Then $|V(Q)|\geq2|Z_3|+2r-5$.\\

\noindent\textbf{Lemma 7}. Let $Q$ be a path and $(Z_{1},Z_{2},Z_{3},Z_{4})$ be a nontrivial $(Q,r)$-scheme with $|Z_{1}|=|Z_{2}|=1$. Then

$$
|V(Q)|\geq\min\bigg\lbrace 2\sum^{4}_{i=1}|Z_{i}|+3r-17,\frac{1}{2}r\Big(\sum^{4}_{i=1}|Z_{i}|-4\Big)+1\bigg\rbrace.
$$

\noindent\textbf{Lemma 8}. Let $Q$ be a path and $(Z_{1},Z_{2},Z_{3},Z_{4})$ be a nontrivial $(Q,r)$-scheme with $|Z_{1}|=|Z_{2}|=|Z_3|=1$ and $|Z_4|\geq4$. Then $|V(Q)|\geq2|Z_4|+3r-7$.\\

Using Woodall's proof technique [15] known as "hopping", we obtain the next result concerning cycles through specified edges.\\

\noindent\textbf{Lemma 9}. Let $G$ be a graph, $A^{\uparrow}$ be a fragment of $G$ with respect to a minimum cut-set $S$ and the connectivity $\kappa$ is even. Let $L$ be a set of $\kappa/2$ independent (vertex disjoint) edges in $\langle S\rangle$ and let $v_1v_2v_3v_4$ be a path in $G$ with $v_1,v_4\in A^{\downarrow}$ and $v_2,v_3\in S$. If a subgraph $\langle S\cup A^{\downarrow}\rangle-\{ v_2,v_3,v_4\}$ contains a cycle $C$ that uses all the edges in $L-\lbrace v_2v_3\rbrace$, then $\langle S\cup A^{\downarrow}\rangle$ contains a cycle that uses all the edges in $L$.\\

In [12, Theorem 1], Veldman proved the following. \\

\noindent\textbf{Lemma B [12]}. If $G$ is a graph with $\delta >3\kappa /2-1$, then no endfragment of $G$ contains a vertex $v$ with $\kappa (G-v)=\kappa -1$.\\

We shall use Lemmas 9 and $B$ to prove the following useful lemma.\\

\noindent\textbf{Lemma 10}. Let $G$ be a 2-connected graph, $A^{\downarrow}$ be an endfragment of $G$ with respect to a minimum cut-set $S$ and let $L$ be a set of independent edges in $\langle S\rangle$. If $\delta > 3\kappa /2 -1$, then $\langle A^{\downarrow}\cup V(L)\rangle$ contains a cycle that uses all the edged in $L$.\\

For the special case $\alpha\leq\delta\leq3\kappa/2-1$ the main lower bound $c\geq\min \{n,4\delta-2\kappa\}$ will be proved by an easy way.\\

\noindent\textbf{Lemma 11}. Every 3-connected graph with $\alpha \leq \delta\leq 3\kappa /2 -1$ has a cycle of length at least $\min \lbrace n,4\delta -2\kappa\}$.\\

Wee need also the following result from [13].\\

\noindent\textbf{Lemma C [13]}. Let $G$ be a hamiltonian graph with $\lbrace v_{1},...,v_{r}\rbrace \subseteq V(G)$ and $d(v_{i})\geq r$ $(i=1,...,r)$. Then any two vertices of $V(G)$ are connected by a path of length at least $r$.\\

Let $V^{\uparrow}$ and $V^{\downarrow}$ are as defined in Definitions A and B. Using above lemmas, we shall prove the following four basic lemmas that are crucial for the proofs of Theorems 2-5.\\

\noindent\textbf{Lemma 12}. Let $G$ be a 3-connected graph with $\delta\geq\alpha$. If $|A^{\uparrow}|\leq 3\delta - \kappa -4$ for a fragment $A^{\uparrow}$ of $G$, then $A^{\uparrow}\subseteq V^{\uparrow}$.\\

\noindent\textbf{Lemma 13}. Let $G$ be a 4-connected graph with $\delta\geq\alpha$. If $|A^{\uparrow}|\geq 3\delta - \kappa -3$ for a fragment $A^{\uparrow}$ of $G$, then either $A^{\uparrow}\subseteq V^{\uparrow}$ or $|V^{\uparrow}|\geq 3\delta -5$.\\

\noindent\textbf{Lemma 14}. Let $G$ be a 3-connected graph with $\delta >3\kappa /2 -1$. If $|A^{\downarrow}|\leq 3\delta - 3\kappa$ for an endfragment $A^{\downarrow}$ of $G$, then $\langle A^{\downarrow}-V^{\downarrow}\rangle$ is edgeless.\\

\noindent\textbf{Lemma 15}. Let $G$ be a 3-connected graph with $\delta >3\kappa /2 -1$ and $|A^{\downarrow}|\geq 3\delta - 3\kappa +1$ for an endfragment $A^{\downarrow}$ of $G$ with respect to a minimum cut-set $S$. If $f=2$ and $S\subseteq V^{\uparrow}$, then either $\langle A^{\downarrow}-V^{\downarrow}\rangle$ is edgeless or $|V^{\downarrow}|\geq 2\delta -2\kappa +3$, where $f=|V^{\downarrow}\cap S|$. If $f=2$ and $S\not\subseteq V^{\uparrow}$, then either $\langle A^{\downarrow}-V(Q^{\downarrow}_0)\rangle$ is edgeless or $|V^{\downarrow}|\geq 3\delta -3\kappa +1$. If $f\geq3$, then either $\langle A^{\downarrow}-V^{\downarrow}\rangle$ is edgeless or $|V^{\downarrow}|\geq 3\delta -3\kappa +f-1$.\\

\section{Proofs of lemmas}

\noindent\textbf{Proof of Lemma 1}. Put $Z=\cup^{3}_{i=1}Z_{i}$. For each $\xi\in Z$, let $f(\xi)$ be the smallest positive integer $h$ such that $\xi^{+h}\in Z$ and let $g(\xi)=|\lbrace i|\xi\in Z_{i}\rbrace |$. Clearly

$$
|C|=\sum_{\xi\in Z}f(\xi),\quad \sum^{3}_{i=1}|Z_{i}|=\sum_{\xi\in Z}g(\xi).\eqno {(1)}
$$

Since $(Z_1,Z_2,Z_3)$ is a nontrivial $(C,r)$-scheme, we have $f(\xi)\geq r$ for each $\xi\in Z$ when $g(\xi)\geq 2$ and $f(\xi)\geq 2$ when $g(\xi)=1$. Let $(\xi_1,\xi_2,\xi_3)$ be a system of distinct representatives of $Z_1,Z_2,Z_3$. Since $|Z_1|=1$, we have $g(\xi)\leq 2$ for each $\xi\in Z-\lbrace \xi_1\rbrace$. In particular, $g(\xi_i)\leq 2$ $(i=2,3)$. Assume first that $r\leq 4$. Clearly $\xi^{+}_{1}\notin Z_1$ and hence $f(\xi_1)\geq r\geq r(g(\xi_1)-1)/2$. Further, for each $\xi\in Z-\lbrace \xi_1\rbrace$, either $g(\xi)=2$ implying that $f(\xi)\geq r=rg(\xi))/2$ or $g(\xi)=1$ implying that $f(\xi)\geq 2\geq rg(\xi)/2$. By summing and using (1), we get 

$$|C|\geq \sum_{\xi\in Z}f(\xi)\geq \Big(\sum^{3}_{i=1}|Z_i|\Big)r/2-r/2$$ 

\noindent and the result follows. Now assume that $r\geq 5$.\\

\textbf{Case 1}. $f(\xi)\geq r$ for each $\xi\in Z$.

By (1), $|C|\geq r|Z|\geq r(|Z_1|+|Z_2|+|Z_3|-1)/2$ and the result follows.\\

\textbf{Case 2}. $f(\xi)\leq r-1$ for some $\xi\in Z$.

Since $g(\xi_2)\leq 2$ and $g(\xi_3)\leq 2$, we can distinguish three subcases.\\

\textbf{Case 2.1}. $g(\xi_2)=g(\xi_3)=1.$

Let $\tau_i$ be the smallest positive integer such that $\xi^{+(\tau_i +1)}_{i}\in Z-Z_i$ and $g(\xi^{+\tau_i}_{i})=1$ $(i=2,3)$. Then\\

$f(\xi^{+\tau_i}_{i})\geq r\geq 2g(\xi^{+\tau_i}_{i})+r-2\quad (i=2,3),$\\

$f(\xi_1)\geq r\geq2g(\xi_1)+r-6.$\\

For each $\xi \in Z-\lbrace \xi^{\tau_2}_{2},\xi^{\tau_3}_{3},\xi_1\rbrace$, either $g(\xi)=2$ implying $f(\xi)\geq r\geq 5 >2g(\xi)$ or $g(\xi)=1$ and again implying $f(\xi)\geq 2=2g(\xi)$. By (1), $|C|\geq 2\sum^{3}_{i=1}|Z_i|+3r-10$ and the result follows.\\

\textbf{Case 2.2}. Either $g(\xi_2)=1$, $g(\xi_3)=2$ or $g(\xi_2)=2$, $g(\xi_3)=1$.

By symmetry, we can assume that $g(\xi_2)=1$ and  $g(\xi_3)=2$. If $g(\xi)=1$ for some $\xi\in Z_3$, then $(\xi_1,\xi_2,\xi)$ is a system of distinct representatives for $(Z_1,Z_2,Z_3)$ and we can argue as in Case 2.1. Let $g(\xi)\geq 2$ for all $\xi\in Z_3$ and let $\tau_2$ be the smallest positive integer such that $\xi^{+(\tau_2 +1)}_{2}\in Z-Z_2$ and $g(\xi^{+\tau_2}_{2})=1$. Then\\

$f(\xi_1)\geq r\geq 2g(\xi_1)+r-6,$\\

$f(\xi^{+\tau_2}_{2})\geq r\geq 2g(\xi^{+\tau_2}_{2})+r-2,\quad f(\xi_3)\geq2g(\xi_3)+r-4.$\\

For each $\xi\in Z-\lbrace \xi_1,\xi^{+\tau_2}_{2},\xi_3\rbrace$, either $g(\xi)=2$ which implies $f(\xi)\geq r\geq 5>2g(\xi)$ or $g(\xi)=1$ again implying $f(\xi)\geq 2=2g(\xi)$. By (1), $|C|\geq 2\sum^{3}_{i=1}|Z_i|+3r-12$ and the result follows.\\

\textbf{Case 2.3.} $g(\xi_2)=g(\xi_3)=2$.

If $g(\xi)=1$ for some $\xi\in Z_2\cup Z_3$, say $\xi\in Z_2 -Z_3$, then $(\xi_1,\xi,\xi_3)$ is a system of distinct representatives and we can argue as in Case 2.2. Otherwise $f(\xi)\geq r$ for all $\xi\in Z$ and we can argue as in Case 1. \quad $\Delta$\\

\noindent\textbf{Proof of Lemma 2}. Put $Z=\cup^{4}_{i=1}Z_i.$ For each $\xi\in Z$, let $f(\xi)$ be the smallest positive integer $h$ such that $\xi^{+h}\in Z$ and let $g(\xi)=|\lbrace i|\xi\in Z_i\rbrace |.$ Then
$$
|C|=\sum_{\xi\in Z}f(\xi),\quad \sum^{4}_{i=1}|Z_i|=\sum_{\xi\in Z}g(\xi). \eqno{(2)}
$$
Let $(\xi_1,\xi_2,\xi_3,\xi_4)$ be a system of distinct representatives of $Z_1,Z_2,Z_3,Z_4$. Since $|Z_1|=|Z_2|=1$, we have $g(\xi)\leq 3$ for each $\xi\in Z$. In particular, $g(\xi_3)\leq 2$ and $g(\xi_4)\leq 2$. Assume first that $r\leq 4$. Clearly $\xi^{+}_{i}\notin Z_i$ $(i=1,2)$ and hence, $f(\xi_i)\geq r\geq r(g(\xi_i)-1)/2$ $(i=1,2)$. For each $\xi\in Z-\lbrace \xi_1,\xi_2\rbrace$, either $g(\xi)=2$ implying $f(\xi)\geq r=rg(\xi)/2$ or $g(\xi)=1$ again implying $f(\xi)\geq 2\geq rg(\xi)/2$. By (2), $|C|\geq r(\sum^{4}_{i=1}|Z_i|)/2-r$ and we are done. Now assume that $r\geq 5$.\\

\textbf{Case 1}. $f(\xi)\geq r$ for each $\xi\in Z$.

By (2), $|C|\geq r|Z|\geq r(\sum^{4}_{i=1}|Z_i|-2)/2$ and the result follows.\\

\textbf{Case 2}. $f(\xi)\leq r-1$ for some $\xi\in Z$.\\

\textbf{Case 2.1.} Either $g(\xi_1)=1$ or $g(\xi_2)=1$.

Assume w.l.o.g. that $g(\xi_1)=1$. Let $p$ be the smallest positive integer such that $\xi^{+p}_{1}\in Z$. Consider two new cycles $C_1$ and $C_2$, obtained from $C$ by identifying $\xi_1$ and $\xi^{+p}_{1}$. Since $f(\xi_1)\geq r$, we have $|C_1|\geq r$ and $|C_2|\leq |C|-r+1.$ Clearly $(Z_2,Z_3,Z_4)$ is a nontrivial $(C_2,r)$-scheme with $|Z_2|=1$. Since $\sum^{4}_{i=1}|Z_i|-1=\sum^{4}_{i=2}|Z_i|$ and $|C|\geq |C_2|+r-1,$ we can obtain the desired result by Lemma 1.\\

\textbf{Case 2.2.} $g(\xi_1)\geq 2$ and $g(\xi_2)\geq 2$.

Clearly $f(\xi_i)\geq r\geq 2g(\xi_i)+r-6$ $(i=1,2)$. If $g(\xi)\geq 2$ for each $\xi\in Z_3 \cup Z_4$, then $f(\xi)\geq r$ for each $\xi\in Z$ and we can argue as in Case 1. Let $g(\xi)=1$ for some $\xi\in Z_3\cup Z_4$. Assume w.l.o.g. that $g(\xi_3)=1$.\\

\textbf{Case 2.2.1.} $g(\xi_4)=1$.

Let $\tau_i$ be  the smallest positive integer such that $\xi^{+(\tau_i +1)}_{i}\in Z-Z_i$ and $g(\xi^{+\tau_i}_{i})=1$ $(i=3,4)$. Then \\

$f(\xi^{+\tau_i}_{i})\geq r\geq 2g(\xi^{+\tau_i}_{i})+r-2\quad (i=3,4).$\\

For each $\xi\in Z-\lbrace \xi^{\tau_3}_{3},\xi^{\tau_4}_{4},\xi_1,\xi_2 \rbrace$, either $g(\xi)=2$ implying $f(\xi)\geq r\geq 5>2g(\xi)$ or $g(\xi)=1$ again implying $f(\xi)\geq 2=2g(\xi)$. By (2), $|C|\geq 2\sum^{4}_{i=1}|Z_i|+4r-16$ and the result follows.\\

\textbf{Case 2.2.2.} $g(\xi_4)=2$.

Let $\tau_3$ be  the smallest positive integer such that $\xi^{+(\tau_3 +1)}_{3}\in Z-Z_i$ and $g(\xi^{+\tau_3}_{3})=1$. Then $f(\xi^{+\tau_3}_{3})\geq r\geq 2g(\xi^{+\tau_3}_{3})+r-2$. If $g(\xi)=1$ for some $\xi\in Z_4$, then we can argue as in Case 2.2.1. Otherwise $g(\xi)=2$ and $f(\xi)\geq r\geq 2g(\xi)+r-4$ for each $\xi\in Z_4 -\lbrace \xi_1,\xi_2\rbrace$. By (2), $|C|\geq 2\sum^{4}_{i=1}|Z_i|+3r-18$ and the result follows. \quad $\Delta$\\

\noindent\textbf{Proofs of Lemmas 3-8}. To prove Lemma 3, form a cycle $C$ consisting of $Q$ and an arbitrary path of length $r$ having only $F(Q)$ and $L(Q)$ in common with $Q$. Since $(Z_1,Z_2)$ is a nontrivial $(C,r)$-scheme, the desired result follows from Lemma A immediately. Lemmas 5 and 7 can be proved by a similar way using Lemmas 1 and 2, respectively. The proofs of Lemmas 4, 6 and 8 are straightforward.\quad $\Delta$\\

\noindent\textbf{Proof of Lemma 9.} We use a variant of an important proof technique known as "hopping" [15]. For the case $v_1\notin V(C)$, we can argue exactly as in [15, proof of Theorem 2]. Let $v_1\in V(C)$. Put $G^{*}=G-\lbrace v_2,v_3\rbrace$ and $L^{\prime}=L-\lbrace v_2v_3\rbrace$. If $X\subseteq V(C)$, we consider all maximal segments of $C-L^{\prime}$ connecting two vertices of $X$. Following [4], the union of the vertex sets of these segments is denoted $Cl(X)$, the endvertices of the segments constitute $Fr(X)$ and finally $Int(X)=Cl(X)-Fr(X)$. The sequence $A_{-1}\subseteq A_0 \subseteq A_{1}\subseteq ...$ of subsets of $V(C)$ is defined as follows: $A_{-1}=\emptyset$ and $A_0$ is the set of vertices $z$ of $C$ such that $G^*$ has a path from $v_4$ to $z$ having only $z$ in common with $C$. For each $p\geq 1$, $A_p$ is the union of $A_{p-1}$ and the set of vertices $z$ such that $G^*$ contains a path $P$ from $Int(A_{p-1})$ to $z$ having only its ends in common with $C$. Let $A=\cup^{\infty}_{i=0}A_i$ and $B=\lbrace v_1\rbrace$. Consider the following statement:

$X(P):$ There exists a path $R_p$ in $G^{*}-\lbrace v_4\rbrace$ starting at $a_p$ in $A_p$ and terminating at $v_1$ such that conditions $(a)-(c)$ below are satisfied.

$(a)$ $R_p$ contains all the edges of $L^{\prime}$ and all the vertices of $Int(A_{p-1})$.

$(b)$ If $Q$ is a segment of $R_p$ from $u$ to $v$ say, having precisely $u$ and $v$ in common with $C$, then one of $u$ and $v$ is outside $A_p$ and the other is outside $\lbrace v_1\rbrace$.

$(c)$ If $y\in Int(X)\cap R_p$, where $X=A_{p^{\prime}}$, $p^{\prime}\leq p-1$ and $M$ denotes the segment of $C-L^{\prime}$ which starts and terminates at $Fr(X)$ and contains $y$ , then $R_p$ contains $M$.

Prove that $X(P)$ holds for some $p$. For suppose this is not the case. Then none of the $\kappa/2-1$ paths of $C-L^{\prime}$ intersects both $A$ and $B$ unless it contains precisely one vertex from $A\cup B$, i.e., $|Fr(A)|\leq 2(\kappa/2-2)+1=\kappa-3.$ Then choose a vertex $z$ on $C$ which is incident with $L^{\prime}$ and not in $Cl(A).$ Now every path in $G^*$ from $v_4$ to $z$ intersects $Fr(A)$, i.e. $|Fr(A)|\geq \kappa -2$. This contradiction proves that $X(P)$ holds for some $p$. Choosing $p$ such that $X(P)$ holds and such that $p$ is minimum, it can be shown that $p=0$ (by the same arguments as in [15, proof of Theorem 2]). Then the desired result holds immediately.\quad $\Delta$\\

\noindent\textbf{Proof of Lemma 10}. The proof is by induction on $\kappa$. For $\kappa =2$ the result follows easily. Let $\kappa \geq 3$. Suppose first that $S-V(L)\neq\emptyset$ and choose any $u\in S-V(L)$. Clearly $A^{\downarrow}$ is an endfragment for $G-u$ too  with respect to $S-u$ and $\delta (G-u)\geq \delta -1> 3\kappa (G-u)/2-1$. By the induction hypothesis, $G-u$ (as well as $G$) contains the desired cycle. Now let $S-V(L)=\emptyset$, i.e. $|L|=k/2$, and choose any $vw\in L$. It follows from $\delta >3\kappa /2-1$ that $|A^{\downarrow}|\geq 2.$ Further, it is not hard to see that there exist two edges $vv^{\prime}$ and $ww^{\prime}$ such that $v^{\prime},w^{\prime}\in A^{\downarrow}$ and $v^{\prime}\neq w^{\prime}$. Put $G^* =G-\lbrace v,w,w^{\prime}\rbrace$ and $S^{\prime}=S-\lbrace v,w\rbrace$. By Lemma B, $\kappa (G-w^{\prime})=\kappa$, i.e. $\kappa (G^*)=\kappa -2$. Also, $\delta (G^*)\geq \delta -3>3\kappa (G^*)/2-1$. If $A^{\downarrow}-\lbrace w^{\prime}\rbrace$ is an endfragment of $G^*$ (with respect to $S^{\prime}$), then by the induction hypothesis, $\langle (S^{\prime}\cup A^{\downarrow})-\lbrace w^{\prime}\rbrace\rangle$ contains a cycle that uses all the edges in $L-\lbrace vw\rbrace$ and the result follows from Lemma 9 immediately. Otherwise choose an endfragment $A^{\downarrow}_{0}\subset A^{\downarrow}-\lbrace w^{\prime}\rbrace$ in $G^*$ with respect to a minimum cut-set $S^{\prime\prime}$ of order $\kappa-2$. Let $P_1,...,P_{\kappa -2}$ be the vertex disjoint paths connecting $S^{\prime}$ and $S^{\prime\prime}$, where $|V(P_i)|=1$ if and only if $F(P_i)=L(P_i)\in S^{\prime}\cap S^{\prime\prime}$ $(i=1,...,\kappa -2)$. By the induction hypothesis, $\langle A^{\downarrow}_{0}\cup S^{\prime\prime}\rangle$ contains a cycle that uses all the independent edges in $\langle S^{\prime\prime}\rangle$ chosen beforehand. Then using $P_1,...,P_{\kappa -2}$, we can form a cycle in $\langle S^{\prime}\cup A^{\downarrow}-\lbrace w^{\prime}\rbrace\rangle$ that uses all the edges in $L-\lbrace vw\rbrace$ and the result follows from Lemma 9. \quad $\Delta$\\

\noindent\textbf{Proof of Lemma 11}. By Theorem G, $c\geq\min\lbrace n,3\delta -3\rbrace$. If $c=n$, then we are done. So, assume that $c\geq 3\delta -3$. If $\kappa \geq 4$, then $c\geq 3\delta -3\geq 4\delta -3\kappa /2-2 \geq 4\delta -2\kappa$. Finally, if $\kappa=3$, then from $\delta\leq3\kappa/2-1$ we get $\delta=3$ implying that $c\geq3\delta-3=4\delta-2\kappa$. \quad $\Delta$\\

\noindent\textbf{Proof of Lemma 12}. Let $S, Q^{\uparrow}_{i}, V^{\uparrow}_{i}$ $(i=1,...,m)$, $Q^{\uparrow}$ and $V^{\uparrow}$ are as defined in Definition A. Assume w.l.o.g. that 

$$
F(Q^{\uparrow}_{i})=u_i, \quad L(Q^{\uparrow}_{i})=v_i \quad (i=1,...,m),
$$

$$
Q^{\uparrow}=u_1\overrightarrow{Q}^{\uparrow}_{1}v_1u_2\overrightarrow{Q}^{\uparrow}_{2}v_2u_3...v_{m-1}u_m\overrightarrow{Q}^{\uparrow}_{m}v_m,
$$

\noindent where $v_1u_2,v_2u_3,...,v_{m-1}u_m$ are extra edges in $G$. Assume the converse, that is, $A^{\uparrow}\not\subseteq V^{\uparrow}$. Let $P=y_1 y_2 ...y_p$ be a longest path in $\langle A^{\uparrow}-V^{\uparrow}\rangle$. Set $Z_1 =N(y_1)\cap V^{\uparrow}$ and $Z_2 =N(y_p)\cap V^{\uparrow}$. Clearly $p+|V^{\uparrow}|\leq |A^{\uparrow}|+|S|\leq 3\delta -4$ and hence

$$
|V^{\uparrow}|\leq 3\delta -p-4. \eqno{(3)}
$$

\textbf{Case 1}. Every path between $V(P)$ and $S-V^{\uparrow}$, intersects $V^{\uparrow}$. 

\textbf{Case 1.1}. $p=1$. 

In this case, $N(y_1)\subseteq V^{\uparrow}$. Set $M=\{u_1,...,u_m\}\cup\{v_1,...,v_m\}$ and $M^*=M\cap N(y_1)$. Since $Q^{\uparrow}_{1},...,Q^{\uparrow}_{m}$ is extreme, $|M^*|\leq 2$. Moreover, $|M^*|=2$ if and only if $M^* =\lbrace u_i,v_i\rbrace$ for some $i \in \lbrace 1,...,m \rbrace$. If $|M^*|\leq 1$, then by standard arguments either $N(y_1)^+ \cup \lbrace y_1 \rbrace$ or $N(y_1)^- \cup \lbrace y_1 \rbrace$ is an independent set of order at lest $\delta+1$, contradicting $\delta\geq \alpha$. So, $|M^*|=2$. Assume w.l.o.g. that $M^*=\{u_1,v_1\}$, that is, $y_1$ is adjacent to both $u_1$ and $v_1$. Put $B=N(y_1)-v_1$. Since $y_1v_i\notin E(G)$ $(i=2,...,m)$, we have $|B^+|\geq\delta-1$. If $B^+\cap S=\emptyset$, then using standard arguments we can show that $B\cup\{y_1,w\}$ for each $w\in A^{\downarrow}$ is an independent set of order at least $\delta+1$, contrary to $\delta\geq\alpha$. Hence, we can choose any $z\in B^+\cap S$. If $z\in V^{\uparrow}_1$, then the collection of paths obtained from $Q^{\uparrow}_{1},...,Q^{\uparrow}_{m}$ by deleting $Q^{\uparrow}_1$ and adding $u_1\overrightarrow{Q}^{\uparrow}_{1}z^-y_1v_1\overleftarrow{Q}^{\uparrow}_{1}z$, contradicts the definition of $Q^{\uparrow}_{1},...,Q^{\uparrow}_{m}$. Therefore, $z\notin V^{\uparrow}_1$. Assume w.l.o.g. that $z\in V^{\uparrow}_2$. If $z\neq v_2$, then we get a new collection of paths obtained from $Q^{\uparrow}_{1},...,Q^{\uparrow}_{m}$ by deleting $Q^{\uparrow}_1$ and $Q^{\uparrow}_2$ and adding $u_1\overrightarrow{Q}^{\uparrow}_{1}v_1y_1z^-\overleftarrow{Q}^{\uparrow}_{2}u_2$ and $z\overrightarrow{Q}^{\uparrow}_{2}v_2$, contrary to $Q^{\uparrow}_{1},...,Q^{\uparrow}_{m}$. So, Let $z=v_2$ implying that $v^-_2\in N(y_1)$. Taking the reverse direction on $Q^{\uparrow}_2$, we can state in addition that $u^+_2\in N(y_1)$. By standard arguments, $B^+\cup\{y_1\}$ is an independent set of vertices of order at least $\delta$. Now we claim that $u_2$ has no neighbors in $B^+\cup\{y_1\}$. Assume, to the contrary, that is, $u_2w\in E(G)$ for some $w\in B^+\cup \{y_1\}$. If $w=y_1$, then deleting $Q^{\uparrow}_1$ and $Q^{\uparrow}_2$ from $Q^{\uparrow}_{1},...,Q^{\uparrow}_{m}$ and adding $u_1\overrightarrow{Q}^{\uparrow}_{1}v_1y_1u_2\overrightarrow{Q}^{\uparrow}_{2}v_2$ we obtain a new collection of paths, contrary to $Q^{\uparrow}_{1},...,Q^{\uparrow}_{m}$. Next, if $w\in V^{\uparrow}_1$, then deleting $Q^{\uparrow}_1$ and $Q^{\uparrow}_2$ and adding $u_1\overrightarrow{Q}^{\uparrow}_{1}w^-y_1u^+_2\overrightarrow{Q}^{\uparrow}_{2}v_2$ and $u_2w\overrightarrow{Q}^{\uparrow}_{1}v_1$ we obtain a new collection, contrary to $Q^{\uparrow}_{1},...,Q^{\uparrow}_{m}$. Further, if $w\in V^{\uparrow}_2$, then deleting $Q^{\uparrow}_1$ and $Q^{\uparrow}_2$ and adding $u_1\overrightarrow{Q}^{\uparrow}_{1}v_1y_1w^-\overleftarrow{Q}^{\uparrow}_{2}u_2w\overrightarrow{Q}^{\uparrow}_{2}v_2$ we obtain a collection, contrary to $Q^{\uparrow}_{1},...,Q^{\uparrow}_{m}$. Finally, if $w\in V^{\uparrow}_i$ for some $i\geq3$, say $i=3$, then deleting $Q^{\uparrow}_2$ and $Q^{\uparrow}_3$ and adding $u_2w\overrightarrow{Q}^{\uparrow}_{3}v_3$ and $u_3\overrightarrow{Q}^{\uparrow}_{3}w^-y_1u^+_2\overrightarrow{Q}^{\uparrow}_{2}v_2$, we obtain a collection, contrary to $Q^{\uparrow}_{1},...,Q^{\uparrow}_{m}$. So, $B^+\cup\{y_1,u_2\}$ is an independent set of order at least $\delta +1$, contrary to $\delta\geq\alpha$.\\

\textbf{Case 1.2}. $p\geq2$. 

If $p=2$, then $|Z_1|\geq \delta -1$, $|Z_2|\geq \delta -1$ and $(Z_1,Z_2)$ is a nontrivial $(Q^{\uparrow},3)$-scheme. By Lemma 3, $|V^{\uparrow}|\geq \min\lbrace 4\delta -6,3\delta -5\rbrace =3\delta -5$, contradicting (3). So, we can assume that $p\geq 3$. Let $w_1,w_2,...,w_s$ be the elements of $(N(y_p)\cap V(P))^+$ occuring on $\overrightarrow{P}$ in a consecutive order, where $w_s=y_p$. Put $P_0=w^{-}_{1}\overrightarrow{P}w_s $ and $p_0=|P_0|$. For each $w_i\in V(P)$ $(i\in \lbrace 1,...,s\rbrace)$ there is a path $y_1 \overrightarrow{P}w^{-}_{i}w_s \overleftarrow{P}w_i$ in $\langle V(P)\rangle$ of length $p$ connecting $y_1$ and $w_i$. Hence, we can assume w.l.o.g. that $P$ is chosen such that for each $i\in\{1,...,s\}$,
 
$$
|Z_1|\geq |N(w_i)\cap V^{\uparrow}|,\quad N(w_i)\cap V(P)\subseteq V(P_0). \eqno{(4)}
$$ 

In particular, $|Z_1|\geq|Z_2|$. Clearly $p_0\geq 2$. If $p_0=2$, then $|Z_1|\geq |Z_2|\geq\delta -1$, and we can argue as in case $p=2$. Let $p_0\geq 3$. Since $G$ is 3-connected, there are vertex disjoint paths $R_1,R_2,R_3$ connecting $P_0$ and $V^{\uparrow}$. Let $F(R_i)\in V(P_0)$ and $L(R_i)\in V^{\uparrow}$ ($i=1,2,3)$.\\ 

\textbf{Case 1.2.1.} $|Z_1|\leq 3.$

By (4), $|N(w_i)\cap V^{\uparrow}|\leq |Z_1|\leq 3$ and therefore, $|N(w_i)\cap V(P_0)|\geq \delta -3$ ($i=1,...,s)$. In particular, for $i=s$, we have $s=|N(y_p)\cap V(P_0)|\geq \delta -3$, implying that $p\geq \delta -2$. By (3), $|V^{\uparrow}|\leq 3\delta -p-4=2\delta-2$. Furthermore, by Lemma C, in $\langle V(P_0)\rangle$ any two vertices are joined by a path  of length at least $\delta -3$.  Due to $R_1,R_2,R_3$, we have $|V^{\uparrow}|\geq 2\delta -1$, a contradiction.\\

\textbf{Case 1.2.2.} $|Z_1|\geq 4$.

Choose $w\in\lbrace w_1,...,w_s\rbrace$ as to maximize $|N(w_i)\cap V^{\uparrow}|$, $i=1,...,s$. Set $Z_3 =N(w)\cap V^{\uparrow}$. By (4), $|Z_1|\geq|Z_3|\geq|Z_2|$. If $|Z_3|\leq3$, then we can argue as in Case 1.2.1. So, we can assume that $|Z_3|\geq 4$. Clearly $|N(w_i)\cap V(P_0)|\geq \delta-|Z_3|$ for each $i\in\{1,...,s\}$. In particular, $s=|N(w_s)\cap V(P_0)|\geq \delta -|Z_3|$ implying that $p\geq \delta-|Z_3|+1$. By Lemma C, in $\langle V(P_0)\rangle$ any two vertices are joined by a path of length at least $\delta-|Z_3|$.\\

\textbf{Case 1.2.1.} $\delta-|Z_3|\geq 1$.

Assume w.l.o.g. that $w\notin V(R_1\cup R_2)$. If $R_1\cup R_2$ does not intersect $y_1 Pw^{-}_{1}$, then $(Z_1,\{L(R_1)\},Z_3)$ is a nontrivial $(Q^{\uparrow}, \delta-|Z_3|+2)$-scheme. Otherwise, let $t$ be the smallest integer such that $y_t\in V(R_1\cup R_2)$. Assume w.l.o.g. that $y_t\in V(R_2)$. Then due to $y_1 \overrightarrow{P}y_tR_2 F(R_2)$, we again can state that $(Z_1,\lbrace L(R_1)\rbrace,Z_3)$ is a nontrivial $(Q^{\uparrow},\delta-|Z_3|+2)$-scheme. By Lemma 5,
$$
|V^{\uparrow}|\geq 2\delta+|Z_3|-4+\min\lbrace |Z_3|-1,(\delta -|Z_3|-1)(|Z_3|-3)\rbrace\geq 2\delta+|Z_3|-4.
$$
On the other hand, using (3) and the fact that $p\geq \delta-|Z_3|+1$, we get $|V^{\uparrow}|\leq 2\delta+|Z_3|-5$, a contradiction.\\

\textbf{Case 1.2.2.2.} $\delta-|Z_3|\leq 0$.

In this case, $|Z_1|\geq|Z_3|\geq\delta$ and $(Z_1,Z_3)$ is a nontrivial $(Q^{\uparrow},p+1)$-scheme. By Lemma 3, $|V^{\uparrow}|\geq3\delta-p-3$, contradicting (3).\\

\textbf{Case 2}. There is a path between $V(P)$ and $S-V^{\uparrow}$ avoiding $V^{\uparrow}$.

Since $Q^{\uparrow}_{1},...,Q^{\uparrow}_{m}$ is extreme, all the paths connecting $V(P)$ and $S-V^{\uparrow}$ and not intersecting $V^{\uparrow}$, end in a unique vertex $z\in S-V^{\uparrow}$.\\
 
\textbf{Case 2.1}. $p=1$.

In this case, $y_1z\in E(G)$ and $N(y_1)-z\subseteq V^{\uparrow}$. Put $B=(N(y_1)-z)^+\cup\{y_1\}$. By standard arguments, $B$ is an independent set of order at least $\delta$. Now we claim that $u_1$ has no neighbors in $B$. Assume, to the contrary, that is, $u_1w\in E(G)$ for some $w\in B$. First, if $w=y_1$, then deleting $Q^{\uparrow}_1$ from $Q^{\uparrow}_1,...,Q^{\uparrow}_m$ and adding $zy_1u_1\overrightarrow{Q}^{\uparrow}_1v_1$ we obtain a new collection of paths, contrary to $Q^{\uparrow}_1,...,Q^{\uparrow}_m$. Next, if $w\in V^{\uparrow}_1$, then deleting $Q^{\uparrow}_1$ and adding $v_1\overleftarrow{Q}^{\uparrow}_1wu_1\overrightarrow{Q}^{\uparrow}_1w^-y_1z$ we obtain another collection of paths, contrary to $Q^{\uparrow}_1,...,Q^{\uparrow}_m$. Finally, if $w\in V^{\uparrow}_i$ for some $i\geq2$, say $i=2$, then deleting $Q^{\uparrow}_1$ and $Q^{\uparrow}_2$ and adding $v_2\overleftarrow{Q}^{\uparrow}_2wu_1\overrightarrow{Q}^{\uparrow}_1v_1$ and $u_2\overrightarrow{Q}^{\uparrow}_2w^-y_1z$ we again obtain a collection, contrary to $Q^{\uparrow}_1,...,Q^{\uparrow}_m$. So, $B\cup \{u_1\}$ is an independent set of order at least $\delta+1$, contrary to $\delta\geq\alpha$.\\

\textbf{Case 2.2}. $p\geq2$.
 
Divide $Q^{\uparrow}$ into three consecutive segments $I_1=\xi_1Q^{\uparrow}\xi_2$, $I_2=\xi_2Q^{\uparrow}\xi_3$ and $I_3=\xi_3Q^{\uparrow}\xi_4$ such that $I_2$ contains $Z_1\cup Z_2$ and is as small as possible. Denote by $R_1$ ($R_2$, respectively) a longest path joining $\xi_2$ ($\xi_3$, respectively) to $z$ and passing through $A^{\uparrow}-V^{\uparrow}$. Since $Q^{\uparrow}_{1},...,Q^{\uparrow}_{m}$ is extreme, $|I_1|\geq|R_1|$ and $|I_2|\geq|R_2|$. Further, we can first estimate $|I_2|$ by Lemma 3 as in Case 1.2, and observing that $|Q^{\uparrow}|=|I_1|+|I_2|+|I_3|\geq|I_2|+|R_1|+|R_2|$, we can argue exactly as in Case 1.2. Lemma 5 can be applied by a similar way with respect to $\cup^{3}_{i=1}Z_i$. \quad $\Delta$\\

\noindent\textbf{Proof of Lemma 13.} Let $S, Q^{\uparrow}_{i}, V^{\uparrow}_{i} (i=1,...,m),$ $Q^{\uparrow}$ and $V^{\uparrow}$ are as defined in Definition A. In addition, let $P=y_1y_2...y_p,P_0,Z_1,Z_2,Z_3$ are as defined in Lemma 12.\\

\textbf{Case 1}. Every path between $V(P)$ and $S-V^{\uparrow}$, intersects $V^{\uparrow}$. 

If $p=0$, then $A^{\uparrow}\subseteq V^{\uparrow}$ and we are done. If $p=1$, then $\alpha\geq\delta+1$ (see the proof of Lemma 12, Case 1.1) contradicting the hypothesis . Further, if $p=2$, then $(Z_1,Z_2)$ is a nontrivial $(Q^{\uparrow},3)$-scheme with $|Z_1|\geq\delta-1,$ $|Z_2|\geq\delta-1$ and, by Lemma 3, $|V^{\uparrow}|\geq3\delta-5$. Now let $p=3$. If $y_1y_3\notin E$ then $|Z_1|\geq\delta-1$, $|Z_2|\geq\delta -1$ and as above, $|V^{\uparrow}|\geq3\delta-5$. Let $y_1 y_3\in E$. This means that $|Z_1|\geq\delta-2$, $|Z_2|\geq\delta-2$ and $(Z_1,Z_2)$ is a nontrivial $(Q^{\uparrow},4)$-scheme. By Lemma 3, $|V^{\uparrow}|\geq4\delta-11$. For $\delta\geq 6$ the inequality $|V^{\uparrow}|\geq3\delta-5$ holds immediately. Let $4\leq\delta\leq5$. Since $\kappa\geq4$, there are four paths connecting $P$ and $Q^{\uparrow}$ and three of them are pairwise disjoint. Then it is easy to see that $|V^{\uparrow}|\geq10\geq3\delta-5$. So, we can assume that $p\geq4$. Suppose first that $|V(P_0)|\leq3$ whence $|Z_1|\geq|Z_2|\geq\delta-2$. Clearly $(Z_1,Z_2)$ is a nontrivial $(Q^{\uparrow},5)$-scheme and by Lemma 3, $|V^{\uparrow}|\geq\min\lbrace 4\delta-8,5\delta-14\rbrace$. If $\delta\geq5$ then $|V^{\uparrow}|\geq3\delta-5$ holds immediately. Otherwise, using 4-connectedness of $G$, it is easy to see that $|V^{\uparrow}|\geq7=3\delta-5$. So, assume that $|V(P_0)|\geq4$. Then $P$ and $Q^{\uparrow}$ are connected by at least four pairwise disjoint paths $R_1,R_2,R_3,R_4$. If $|Z_1|\leq3$, then as in Lemma 12 (Case 1.2.1), in $\langle V(P_0)\rangle$ each two vertices are connected by a path of length at least $\delta-3$ and due to $R_1,R_2,R_3,R_4$, $|V^{\uparrow}|\geq3(\delta-1)+1>3\delta-5$. Let $|Z_1|\geq4$. By similar arguments, $|Z_3|\geq4$. Clearly $|N(w_i)\cap V(P_0)|\geq\delta-|Z_3|$ $(i=1,...,s)$ and by Lemma C, in $\langle V(P_0)\rangle$ any two vertices are joined by a path of length at least $\delta-|Z_3|$. If $\delta-|Z_3|\geq1$, then we can assume w.l.o.g. (see the proof of Lemma 12, Case 1.2.2.1) that $(\lbrace L(R_1)\rbrace,\lbrace L(R_2)\rbrace,Z_1,Z_3)$ is a nontrivial $(Q^{\uparrow},\delta-|Z_3|+2)$-scheme and by Lemma 7,\\

$|V^{\uparrow}|\geq3\delta-5+\min\lbrace |Z_3|-4,(\delta-|Z_3|-1)(|Z_3|-4)\rbrace\geq3\delta-5.$\\

Otherwise $|Z_1|\geq|Z_3|\geq\delta$ and $(Z_1,Z_3)$  is a nontrivial $(Q^{\uparrow},p+1)$-scheme. By Lemma 3,\\

$|V^{\uparrow}|\geq3\delta-5+\min\lbrace\delta+p+1,(\delta-1)(p-2)+3\rbrace\geq3\delta-5.$\\

\textbf{Case 2}. There is a path between $V(P)$ and $S-V^{\uparrow}$ avoiding $V^{\uparrow}$.

We can argue exactly as in proof of Lemma 12 (Case 2). \quad $\Delta$\\

\noindent\textbf{Proof of Lemma 14}. Let $S, Q^{\downarrow}_{1},...,Q^{\downarrow}_{m}$ and $V^{\downarrow}_{1},...,V^{\downarrow}_{m},V^{\downarrow}$ are as defined in Definition B. The existence of $Q^{\downarrow}_{1},...,Q^{\downarrow}_{m}$ follows from Lemma 10. Put $|V^{\downarrow}\cap S|=f$. Clearly $f\geq2m$. We can assume that $\delta-\kappa\geq2$ since otherwise $|A^{\downarrow}|\leq3$ (by the hypothesis) and it is not hard to see that $\langle A^{\downarrow}-V^{\downarrow}\rangle$ is edgeless. Let $P=y_1y_2...y_p$ be a longest path in $\langle A^{\downarrow}-V^{\downarrow}\rangle$. By the hypothesis, $p+|V^{\downarrow}|-f\leq |A^{\downarrow}|\leq3\delta-3\kappa$ implying that

$$|V^{\downarrow}|\leq3\delta-3\kappa-p+f. \eqno{(5)}$$

Put\\

$Z_1=N(y_1)\cap V^{\downarrow},\quad Z_2=N(y_p)\cap V^{\downarrow},$\\

$Z_{1,i}=Z_1\cap V^{\downarrow}_{i},\quad Z_{2,i}=Z_2\cap V^{\downarrow}_{i}\quad (i=1,...,m).$\\

Clearly $Z_1=\cup^{m}_{i=1}Z_{1,i}$ and $Z_2=\cup^{m}_{i=1}Z_{2,i}$. If $p\leq1$, then $\langle A^{\downarrow}-V^{\downarrow}\rangle$ is edgeless and we are done. Let $p\geq2$.\\

\textbf{Case 1}. $p=2$.

In this case, $|Z_i|\geq\delta-\kappa+f-1$ $(i=1,2)$. We claim that 

$$
|V^{\downarrow}_{i}|\geq\frac{3}{2}(|Z_{1,i}|+|Z_{2,i}|)-2 \quad (i=1,...,m). \eqno{(6)}
$$

Indeed, if $(Z_{1i},Z_{2i})$ is a nontrivial $(Q^{\downarrow}_{i},3)$-scheme, then (6) holds by Lemma 3, immediately. Otherwise it can be checked easily. By summing,
$$
|V^{\downarrow}|=\sum^{m}_{i=1}|V^{\downarrow}_{i}|\geq\frac{3}{2}(|Z_1|+|Z_2|)-2m\geq3\delta-3\kappa-p+f+1,
$$
contradicting (5).\\

\textbf{Case 2.} $p\geq3$.

Let $w_1,w_2,...,w_s$ be the elements of $(N(y_p)\cap V(P))^+$ occuring on $\overrightarrow{P}$ in a consecutive order. Put $P_0=w^-_1\overrightarrow{P}w_s$ and $p_0=|V(P_0)|$. As in proof of Lemma 12 (see (4)), we can assume w.l.o.g. that for each $i\in\lbrace1,...,s\rbrace$,
$$
|Z_1|\geq|N(w_i)\cap V^{\downarrow}|,\quad N(w_i)\cap V(P)\subseteq V(P_0).\eqno{(7)}
$$ 
Choose $w\in\lbrace w_1,...,w_s\rbrace$ as to maximize $|N(w_i)\cap V^{\downarrow}|$, $i=1,...,s$. Set
$$
Z_3=N(w)\cap V^{\downarrow}, \quad Z_{3,i}=Z_3\cap V^{\downarrow}_{i} \quad (i=1,...,m).
$$

Clearly $|Z_1|\geq|Z_2|\geq\delta-\kappa+f-p_0+1$ and $|Z_3|\geq|Z_2|\geq\delta-\kappa+f-p_0+1$. We claim that 
 
$$
|V^{\downarrow}_{i}|\geq2(|Z_{1,i}|+|Z_{2,i}|)-3 \quad (i=1,...,m).\eqno{(8)}
$$

Indeed, if $(Z_{1,i},Z_{2,i})$ is a nontrivial $(Q^{\downarrow}_{i},p+1)$-scheme, then (8) holds by Lemma 3 and the fact that $p\geq3$ . Otherwise, it can be checked easily. Analogously, 

$$
|V^{\downarrow}_{i}|\geq2(|Z_{1,i}|+|Z_{3,i}|)-3 \quad (i=1,...,m).\eqno{(9)}
$$

\textbf{Case 2.1.} $p_0\leq m+1$.

Using (8) and summing, we get\\

$|V^{\downarrow}|=\sum^{m}_{i=1}|V^{\downarrow}_{i}|\geq2(|Z_1|+|Z_2|)-3m\geq4(\delta-\kappa+f+p_0+1)-3m$\\

$=(3\delta-3\kappa-p+f+1)+(\delta-\kappa)+(p-p_0)+3(m-p_0+1)+3(f-2m)$.\\

Recalling that $\delta-\kappa\geq2$, $p\geq p_0$, $m-p_0+1\geq0$ and $f\geq2m$, we get $|V^{\downarrow}|\geq3\delta-3\kappa-p+f+1$, contradicting (5).\\

\textbf{Case 2.2.} $p_0 \geq m+2$.

Assume first that $\delta-\kappa+f-|Z_3|\leq1$. Then $|Z_1|\geq|Z_3|\geq\delta-\kappa+f-1$. Using (8) and summing, we get\\ 

$|V^{\downarrow}|=\sum^{m}_{i=1}|V^{\downarrow}_{i}|\geq4(\delta-\kappa+f-1)-3m$\\

$\geq(3\delta-3\kappa-p+f+1)-(\delta-\kappa-2)+3(f-m-1)\geq3\delta-3\kappa-p+f+1$,\\

\noindent contradicting (5). Now assume that $\delta-\kappa+f-|Z_3|\geq2$. Clearly $|N(w_i)\cap V(P_0)|\geq\delta-\kappa+f-|Z_3|$ $(i=1,...,s)$. In particular, for $i=s$, we have $s\geq\delta-\kappa+f-|Z_3|$. By Lemma C, in $\langle V(P_0)\rangle$ any two vertices are joined by a path of length at least $\delta-\kappa+f-|Z_3|$. Observing that $p\geq s+1\geq\delta-\kappa+f-|Z_3|+1$ and combining it with (5), we get

$$
|V^{\downarrow}|\leq2\delta-2\kappa+|Z_3|-1. \eqno {(10)}
$$

\textbf{Case 2.2.1.} $p_0\leq f$.

Since $\kappa\geq f\geq p_0$, there are vertex disjoint paths $R_1,...,R_{p_0}$ connecting $V(P_0)$ and $V^{\downarrow}$. Let $F(R_i)\in V(P_0)$ and $L(R_i)\in V^{\downarrow}$ $(i=1,...,p_0)$. Since $p_0\geq m+2$, we can assume w.l.o.g. that either $L(R_i)\in V^{\downarrow}_{1}$ $(i=1,2)$ and $L(R_i)\in V^{\downarrow}_{2}$ $(i=3,4)$ or $L(R_i)\in V^{\downarrow}_{1}$ $(i=1,2,3)$.\\

\textbf{Case 2.2.1.1}. $L(R_i)\in V^{\downarrow}_{1}$ $(i=1,2)$ and $L(R_i)\in V^{\downarrow}_{2}$ $(i=3,4)$.

Assume w.l.o.g. that $R_1\cup R_3$ does not intersect $y_1\overrightarrow{P}w^{-}_{1}$ (see the proof of Lemma 12, Case 1.2.2.1). Because of the symmetry, we can distinguish the following six subcases.\\

\textbf{Case 2.2.1.1.1}. $|Z_{11}|\leq2$, $|Z_{12}|\leq2$, $|Z_{31}|\leq2$, $|Z_{32}|\leq2$.

Due to $R_1,R_2$ and $R_3,R_4$ we have $|V^{\downarrow}_{i}|\geq\delta-\kappa+f-|Z_3|+3$ $(i=1,2)$. Using (9) for each $i\in\{3,...,m\}$ and summing, we get\\

$|V^{\downarrow}|=|V^{\downarrow}_{1}|+|V^{\downarrow}_{2}|+\sum^{m}_{i=3}|V^{\downarrow}_{i}|$\\

$\geq2(\delta-\kappa+f-|Z_3|+3)+2(|Z_1|-|Z_{11}|-|Z_{12}|+|Z_3|-|Z_{31}|-|Z_{32}|)-3(m-2)$\\ 

$=(2\delta-2\kappa+|Z_1|)+(|Z_1|-|Z_{11}|-|Z_{12}|)+(2f-3m)$\\

$+12-(|Z_{11}|+|Z_{12}|+2|Z_{31}|+2|Z_{32}|)>2\delta-2\kappa+|Z_3|$,\\

\noindent contradicting (10).\\

\textbf{Case 2.2.1.1.2}. $|Z_{11}|\geq3$, $|Z_{12}|\leq2$, $|Z_{31}|\leq2$, $|Z_{32}|\leq2$.

Clearly either $(\{L(R_1)\},Z_{11})$ or $(\{L(R_2)\},Z_{11})$ is a nontrivial $(Q^{\downarrow}_{1},\delta-\kappa+f-|Z_3|+2)$-scheme. By Lemma 4, $|V^{\downarrow}_{1}|\geq(\delta-\kappa +f-|Z_3|+2)+2|Z_{11}|-3$. Due to $R_1$ and $R_2$ we have $|V^{\downarrow}_{2}|\geq\delta-\kappa +f-|Z_3|+3$. Using also (9) for each $i\in\{3,...,m\}$ and summing, we get\\

$|V^{\downarrow}|\geq2(\delta-\kappa+f-|Z_3|+2)+2|Z_{11}|-2+\sum^{m}_{i=3}(2(|Z_{1i}|+|Z_{3i}|)-3)$\\

$=(2\delta-2\kappa+|Z_1|)+(|Z_1|-|Z_{11}|-|Z_{12}|)+2f-3m$\\

$+(|Z_{11}|+8)-(|Z_{12}|+2|Z_{31}|+2|Z_{32}|)>2\delta-2\kappa+|Z_3|$,\\
 
\noindent contradicting (10).\\

\textbf{Case 2.2.1.1.3}. $|Z_{11}|\geq3$, $|Z_{12}|\geq3$, $|Z_{31}|\leq2$, $|Z_{32}|\leq2$.

Clearly either $(\{L(R_1)\},Z_{11})$ or $(\{L(R_2)\},Z_{11})$ is a nontrivial $(Q^{\downarrow}_{1},\delta-\kappa+f-|Z_3|+2)$-scheme. By the same reason, either $(\{L(R_3)\},Z_{12})$ or $(\{L(R_4)\},Z_{12})$ is a nontrivial $(Q^{\downarrow}_{2},\delta-\kappa+f-|Z_3|+2)$-scheme too. By Lemma 4,

$$
|V^{\downarrow}_{i}|\geq(\delta-\kappa+f-|Z_3|+2)+2|Z_{1i}|-3 \quad (i=1,2).
$$

Using (9) for each $i\in\{3,...,m\}$ and summing, \\

$|V^{\downarrow}|\geq|V^{\downarrow}_{1}|+|V^{\downarrow}_{2}|+\sum^{m}_{i=3}|V^{\downarrow}_{i}|\geq2(\delta-\kappa+f-|Z_3|+2)$\\

$+2(|Z_{11}|+|Z_{12}|)-6+2(|Z_1|-|Z_{11}|-|Z_{12}|+|Z_3|-|Z_{31}|-|Z_{32}|)-3(m-2)$\\

$=(2\delta-2\kappa+|Z_1|)+2f-3m+(|Z_1|+4)-2(|Z_{31}|+|Z_{32}|)\geq2\delta-2\kappa+|Z_3|$,\\

\noindent contradicting (10).\\

\textbf{Case 2.2.1.1.4}. $|Z_{11}|\geq3$, $|Z_{12}|\leq2$, $|Z_{31}|\geq3$, $|Z_{32}|\leq2$.

Since $(Z_{11},Z_{31})$ is a nontrivial $(Q^{\downarrow}_{1},\delta-\kappa+f-|Z_3|+2)$-scheme, we can apply Lemma 3,\\

$|V^{\downarrow}_{1}|\geq2(|Z_{11}|+|Z_{31}|)+(\delta-\kappa+f-|Z_3|)-5+$\\

$\min\{4,\frac{1}{2}(\delta-\kappa+f-|Z_3|-2)(|Z_{11}|+|Z_{31}|-6)\}$\\

$\geq(\delta-\kappa+f-|Z_3|)+2(|Z_{11}|+|Z_{31}|)-5.$\\

Due to $R_3$ and $R_4$, we have $|V^{\downarrow}_{2}|\geq\delta-\kappa+f-|Z_3|+3$. Using (9) for each $i\in\{3,...,m\}$ and summing, we get\\

$|V^{\downarrow}|\geq2(\delta-\kappa+f-|Z_3|)+2(|Z_{11}|+|Z_{31}|)-2$\\

$+2(|Z_1|-|Z_{11}|-|Z_{12}|+|Z_3|-|Z_{31}|-|Z_{32}|)-3(m-2)$\\

$=(2\delta-2\kappa+|Z_1|)+2f-3m+(|Z_1|+4)-2(|Z_{12}|+|Z_{32}|)\geq2\delta-2\kappa+|Z_3|,$\\
 
\noindent contradicting (10).\\

\textbf{Case 2.2.1.1.5}. $|Z_{11}|\geq3$, $|Z_{12}|\geq3$, $|Z_{31}|\geq3$, $|Z_{32}|\leq2$.

As in Case 2.2.1.1.4, $|V^{\downarrow}_{1}|\geq(\delta-\kappa+f-|Z_3|)+2(|Z_{11}|+|Z_{31}|)-5$. By Lemma 4, $|V^{\downarrow}_{2}|\geq(\delta-\kappa+f-|Z_3|+2)+2|Z_{12}|-3$. Using (9) for each $i\in\{3,...,m\}$ and summing, we get\\

$|V^{\downarrow}|\geq2(\delta-\kappa+f-|Z_3|)+2(|Z_{11}|+|Z_{31}|)+2|Z_{12}|-6$\\

$+2(|Z_1|-|Z_{11}|-|Z_{12}|+|Z_3|-|Z_{31}|-|Z_{3i}|)-3(m-2)$\\

$\geq(2\delta-2\kappa+|Z_1|)+2f-3m+|Z_1|-|Z_{32}|>2\delta-2\kappa+|Z_3|$,\\

\noindent contradicting (10).\\

\textbf{Case 2.2.1.1.6}. $|Z_{11}|\geq3$, $|Z_{12}|\geq3$, $|Z_{31}|\geq3$, $|Z_{32}|\geq3$.

As in Case 2.2.1.1.4, \\

$|V^{\downarrow}_{1}|+|V^{\downarrow}_{2}|\geq2(\delta-\kappa+f-|Z_3|)+2(|Z_{11}|+|Z_{31}|)+2(|Z_{12}|+|Z_{32}|)-10.$\\

Using (9) for each $i\in\{3,...,m\}$ and summing, we get\\

$|V^{\downarrow}|\geq2(\delta-\kappa+f-|Z_3|)+2(|Z_{11}|+|Z_{31}|+|Z_{12}|+|Z_{32}|)-10$\\

$+2(|Z_1|-|Z_{11}|-|Z_{12}|+|Z_3|-|Z_{31}|-|Z_{32}|)-3(m-2)$\\

$=(2\delta-2\kappa+|Z_1|)+2f-3m+|Z_1|-4\geq2\delta-2\kappa+|Z_3|,$\\

\noindent contrary to (10).\\

\textbf{Case 2.2.1.2}. $L(R_i)\in V^{\downarrow}_{1}$ $(i=1,2,3)$.

Assume w.l.o.g. that $R_1\cup R_2$ does not intersect $y_1\overrightarrow{P}w^{-}_{1}$ (see the proof of Lemma 12, Case 1.2.2.1).\\

\textbf{Case 2.2.1.2.1}. $|Z_{11}|\leq3$, $|Z_{31}|\leq3$.

Due to $R_1,R_2,R_3$ we have $|V^{\downarrow}_{1}|\geq2(\delta-\kappa+f-|Z_3|+2)+1$. Using (9) for each $i\in\{2,...,m\}$ and summing, we get\\

$|V^{\downarrow}|\geq|V^{\downarrow}_{1}|+\sum^{m}_{i=2}|V^{\downarrow}_{i}|$\\

$\geq2(\delta-\kappa+f-|Z_3|+2)+1+2(|Z_1|-|Z_{11}|+|Z_3|-|Z_{31}|)-3(m-1)$\\

$=(2\delta-2\kappa+|Z_1|)+(|Z_1|-|Z_{11}|)+(2f-3m)$\\

$+8-(|Z_{11}|+2|Z_{31}|)\geq2\delta-2\kappa+|Z_3|,$\\

\noindent which contradicts (10).\\

\textbf{Case 2.2.1.2.2}. $|Z_{11}|\leq3$, $|Z_{31}|\geq4$.

Since $|Z_{31}|\geq4$, we can suppose w.l.o.g. that $F(R_3)=w$. Then clearly $(\{L(R_1)\},\{L(R_2)\},Z_{31})$ is a nontrivial $(Q^{\downarrow}_{1},\delta-\kappa+f-|Z_3|+2)$-scheme and by Lemma 6, $|V^{\downarrow}_{1}|\geq2(\delta-\kappa+f-|Z_3|+2)+2|Z_{31}|-5$. Using (9) for each $i\in\{2,...,m\}$ and summing, we get\\

$|V^{\downarrow}|\geq2(\delta-\kappa+f-|Z_3|+2)+2|Z_{31}|-5+2(|Z_1|-|Z_{11}|+|Z_3|-|Z_{31}|)-3(m-1)$\\

$=(2\delta-2\kappa+|Z_1|)+(|Z_1|-|Z_{11}|)+2f-3m+2-|Z_{11}|\geq2\delta-2\kappa+|Z_3|,$\\
 
\noindent contradicting (10).\\

\textbf{Case 2.2.1.2.3}. $|Z_{11}|\geq4$, $|Z_{31}|\geq4$.

For some $R_1,R_2,R_3$, say $R_1$, we have $F(R_1)\notin\{y_1,w\}$. Then clearly $(\{L(R_1)\},Z_{11},Z_{31})$ is a nontrivial $(L^{\downarrow}_{1},\delta-\kappa+f-|Z_3|+2)$-scheme. By Lemma 5,\\

$|V^{\downarrow}_{1}|\geq\min\{2(\delta-\kappa+f-|Z_3|+2)+2(|Z_{11}|+|Z_{31}|+1)-11,$\\

$(\delta-\kappa+f-|Z_3|+2)(|Z_{11}|+|Z_{31}|-2)/2+1\}$\\

$=2(\delta-\kappa+f-|Z_3|)+2(|Z_{11}|+|Z_{31}|)-7$\\

$+\min\{2,(\delta-\kappa+f-|Z_3|-2)(|Z_{11}|+|Z_{31}|-6)/2\}$\\

$\geq2(\delta-\kappa+f-|Z_3|)+2(|Z_{11}|+|Z_{31}|)-7.$\\

Further, applying (9) for each $i\in\{2,...,m\}$ and summing, we get\\

$|V^{\downarrow}|\geq2(\delta-\kappa+f-|Z_3|)+2(|Z_{11}|+|Z_{31}|)-7$\\

$+2(|Z_1|-|Z_{11}|+|Z_3|-|Z_{31}|)-3(m-1)$\\

$=(2\delta-2\kappa+|Z_1|)+|Z_1|-4+2f-3m\geq2\delta-2\kappa+|Z_3|,$\\

\noindent contradicting (10).\\

\textbf{Case 2.2.2}. $p_0\geq f+1$.

Let $S=\{v_1,...,v_{\kappa}\}$ and $V^{\downarrow}\cap S=\{v_1,...,v_{\kappa}\}$. Consider a new graph $G^{\prime}=G-\{v_{f+1},v_{f+2},...,v_\kappa\}$. Add new vertices $a_1,a_2$ in $G^{\prime}$ and join $a_1$ to all vertices of $V^{\downarrow}$, and join $a_2$ to all vertices of $V(P_0)$. Set $G^{\prime\prime}=\langle V(G^{\prime})\cup \{a_1,a_2\}\rangle$. Clearly $G^{\prime\prime}$ is $f$-connected. Let $V_0$ be a minimum cut-set in $G^{\prime\prime}$ that separates $a_1$ and $a_2$. Since $a_1$ and $a_2$ are connected in $G^{\prime\prime}-\{v_1,...,v_f\}$, we have $V_0\neq\{v_1,...,v_f\}$. Observing also that $A^{\downarrow}$ is an endfragment for $G^{\prime}$, we can suppose that $|V_0|\geq f+1$ and therefore there exist $f+1$ internally disjoint paths in $G^{\prime\prime}$ joining $a_1$ and $a_2$. This means that in $G^{\prime}$ (as well as in $G$) there exist vertex disjoint paths $R_1,R_2,...,R_{f+1}$ connecting $V^{\downarrow}$ and $V(P_0)$. Then using the fact that $f+1\geq m+2$, we can argue exactly as in Case 2.2.1.  \quad $\Delta$\\

\noindent\textbf{Proof of Lemma 15}. Let $Q^{\downarrow}_{1},...,Q^{\downarrow}_{m}$ and $V^{\downarrow}_{1},...,V^{\downarrow}_{m}$, $V^{\downarrow}$ be as defined in Definition B. The existence of $Q^{\downarrow}_{1},...,Q^{\downarrow}_{m}$ follows from Lemma 10. Further, let $P=y_1...y_p$, $Z_1$, $Z_2$ and $Z_{1,i}$, $Z_{2,i}$ $(i=1,...,m)$ be as defined in Lemma 14. \\

\textbf{Case 1}. $f\geq3$.

Suppose first that $\delta-\kappa\leq1$. Combining it with $\delta>3\kappa /2-1$, we get $\kappa=3$, $\delta=4$, $f=3$ and $m=1$. Then it is easy to show that $|V^{\downarrow}|\geq5=3\delta-3\kappa+f-1$ and we are done. Now let $\delta-\kappa\geq2$. If $p\leq1$, then clearly $<A^{\downarrow}-V^{\downarrow}>$ is edgeless and we are done. Further, if $p=2$, then $|V^{\downarrow}|\geq3\delta-3\kappa+f-1$ (see the proof of Lemma 14, Case 1). Let $p\geq3$.\\

\textbf{Case 1.1.} $p=3$.

In this case, $P=y_1y_2y_3$. If $y_1y_3\notin E(G)$, then we can argue as in case $p=2$. Let $y_1y_3\in E(G)$ implying that $|Z_i|\geq\delta-\kappa+f-2$ ($i=1,2)$. Applying (8) (see the proof of Lemma 14) and summing, we get \\

$|V^{\downarrow}|\geq\sum^{m}_{i=1}|V^{\downarrow}_{i}|\geq2(|Z_1|+|Z_2|)-3m\geq4(\delta-\kappa+f-2)-3m$\\

$=(3\delta-3\kappa+f-1)+(3f-3m-5)+\delta-\kappa-2.$\\

If $f=3$, then $m=1$ and $3f-3m-5\geq1$. If $f\geq4$, then $3f-3m-5\geq f+m-5\geq0$. In both cases the desired result follows immediately.\\

\textbf{Case 1.2}. $p\geq4$.

Let $P_0,p_0,w,Z_3$ and $Z_{3,i}$ $(i=1,...,m)$ are as defined in Lemma 14 (Case 2). Clearly $|Z_1|\geq|Z_2|\geq\delta-\kappa+f-p_0+1$ and $|Z_1|\geq|Z_3|\geq\delta-\kappa+f-p_0+1$. By the definition of $Q^{\uparrow}_{1},...,Q^{\uparrow}_{m}$, we have $|V^{\downarrow}_{i}|\geq3$ $(i=1,...,m)$. We claim that 

$$
|V^{\downarrow}_{i}|\geq2(|Z_{1,i}|+|Z_{2,i}|)-2 \quad (i=1,...,m). \eqno{(11)} 
$$

Indeed, if $(Z_{1,i},Z_{2,i})$ is not a nontrivial $(Q^{\downarrow}_{i},p+1)$-scheme, then (11) can bee checked easily. Otherwise, by Lemma 3,\\

$|V^{\downarrow}_{i}|\geq\min\{2(|Z_{1,i}|+|Z_{2,i}|)+p-4,\frac{1}{2}(p+1)(|Z_{1,i}|+|Z_{2,i}|-2)+1\}$\\

$=2(|Z_{1,i}|+|Z_{2,i}|)-2+\min\{p-2,\frac{1}{2}(p-3)(|Z_{1,i}|+|Z_{2,i}|-4)+p-4\}.$\\

If $|Z_{1,i}|+|Z_{2,i}|\geq4$, then (11) holds immediately. If $|Z_{1,i}|+|Z_{2,i}|=3$, then (11) follows from $|V^{\downarrow}_{i}|\geq2(|Z_{1,i}|+|Z_{2,i}|)-5/2$. Finally, if $|Z_{1,i}|+|Z_{2,i}|\leq2$, then (11) follows from $|V^{\downarrow}_{i}|\geq3$. By a similar way, 

$$
|V^{\downarrow}_{i}|\geq2(|Z_{1,i}|+|Z_{3,i}|)-2 \quad (i=1,...,m). \eqno{(12)} 
$$

If $p_0\leq m+1$, then using (11) and summing, we get\\ 

$|V^{\downarrow}|=\sum^{m}_{i=1}|V^{\downarrow}_{i}|\geq2(|Z_1|+|Z_2|)-2m\geq4(\delta-\kappa+f-p_0+1)-2m$\\

$=(3\delta-3\kappa+f-1)+\delta-\kappa+3f-4p_0-2m+5$\\

$\geq(3\delta-3\kappa+f-1)+4(m-p_0+1)+(\delta-\kappa+1)>3\delta-3\kappa+f-1.$\\

Now let $p_0\geq m+2$. If $\delta-\kappa+f-|Z_3|\leq1$, then $|Z_1|\geq|Z_3|\geq\delta-\kappa+f-1$ and by (11),\\

$|V^{\downarrow}|=\sum^{m}_{i=1}|V^{\downarrow}_{i}|\geq2(|Z_1|+|Z_2|)-2m\geq4(\delta-\kappa+f-1)-2m$\\

$=(3\delta-3\kappa+f-1)+(\delta-\kappa-2)+3f-2m-2>3\delta-3\kappa+f-1.$\\

Let $\delta-\kappa+f-|Z_3|\geq2$. By the choice of $w$,  \\

$|N(w_i)\cap V(P_0)|\geq\delta-\kappa+f-|Z_3|\quad (i=1,...,s)$. \\

In particular, for $i=s$, we have $s\geq\delta-\kappa+f-|Z_3|$. By Lemma C, in $\langle V(P_0)\rangle$ any two vertices are joined by a path of length at least $\delta-\kappa+f-|Z_3|$.\\

\textbf{Case 1.2.1.} $p_0\leq f$.

Let $G^{\prime}=G-(S-V^{\downarrow})$. Since $G^{\prime}$ is $p_0$-connected, there exist vertex disjoint paths $R_1,R_2,...,R_{p_0}$ connecting $V(P_0)$ and $V^{\downarrow}$. Let $F(R_i)\in V(P_0)$ and $L(R_i)\in V^{\downarrow}$ $(i=1,...,p_0)$. Since $p_0\geq m+2$, we can assume w.l.o.g. that either $L(R_i)\in V^{\downarrow}_{1}$ $(i=1,2)$ and  $L(R_i)\in V^{\downarrow}_{2}$ $(i=3,4)$ or $L(R_i)\in V^{\downarrow}_{1}$ $(i=1,2,3)$.\\

\textbf{Case 1.2.1.1.} $L(R_i)\in V^{\downarrow}_{1}$ $(i=1,2)$ and $L(R_i)\in V^{\downarrow}_{2}$ $(i=3,4)$.

Assume w.l.o.g. that $R_1\cup R_2$ does not intersect $w^{-}_{1}P_0w_s$ (see the proof of Lemma 12, Case 1.2.2.1).\\

\textbf{Case 1.2.1.1.1.} $|Z_{11}|\leq2$, $|Z_{12}|\leq2$, $|Z_{31}|\leq2$, $|Z_{32}|\leq2$.

Due to $R_1,R_2$ and $R_3,R_4$, we have  $|V^{\downarrow}_{i}|\geq\delta-\kappa+f-|Z_3|+3$ $(i=1,2)$. Using (12) for each $i\in\{3,...,m\}$ and summing, we get\\ 

$|V^{\downarrow}|=|V^{\downarrow}_{1}|+|V^{\downarrow}_{2}|+\sum^{m}_{i=3}|V^{\downarrow}_{i}|$\\

$\geq2(\delta-\kappa+f-|Z_3|+3)+2(|Z_1|-|Z_{11}|-|Z_{12}|+|Z_3|-|Z_{31}|-|Z_{32}|)-2(m-2)$\\

$=(3\delta-3\kappa+f-1)+(|Z_1|-|Z_{11}|-|Z_{12}|)+(|Z_3|-\delta+\kappa-f+p_0-1)+(|Z_1|-|Z_3|)$\\

$+(2f-2m-p_0)+12-|Z_{11}|-|Z_{12}|-2|Z_{31}|-2|Z_{32}|$.\\

Since $|Z_3|\geq\delta-\kappa+f-p_0+1$ and $2f-2m-p_0\geq f-p_0\geq0$, we have $|V^{\downarrow}|\geq3\delta-3\kappa+f-1$.\\

\textbf{Case 1.2.1.1.2}. Either $|Z_{11}|\geq3$ or $|Z_{12}|\geq3$ or $|Z_{31}|\geq3$ or $|Z_{32}|\geq3$.

In this case we can argue exactly as in proof of lemma 14 (Cases 2.2.1.1.2-2.2.1.1.6).\\

\textbf{Case 1.2.1.2}. $L(R_i)\in V^{\downarrow}_{1}$ $(i=1,2,3)$.

Assume w.l.o.g. that $R_1\cup R_2$ does not intersect $y_1\overrightarrow{P}w^-_1$ (see the proof of Lemma 12, Case 1.2.2.1).\\

\textbf{Case 1.2.1.2.1}. $|Z_{11}|\leq3$, $|Z_{31}|\leq3$.

Due to $R_1$, $R_2$, $R_3$, we have $|V^{\downarrow}_1|\geq2(\delta-\kappa+f-|Z_3|+2)+1$. Using (12) for each $i\in\{2,...,m\}$ and summing, we get\\

$|V^{\downarrow}|=|V^{\downarrow}_1|+\sum^{m}_{i=2}|V^{\downarrow}_i|$\\

$\geq2(\delta-\kappa+f-|Z_3|+2)+1+2(|Z_1|-|Z_{11}|+|Z_3|-|Z_{31})-2(m-1)$\\

$=(3\delta-3\kappa+f-1)+(|Z_3|-\delta+\kappa-f+p_0-1)$\\

$+(|Z_1|-|Z_3|)+(|Z_1|-|Z_{11}|)+9-|Z_{11}|-2|Z_{31}|+2f-2m-p_0$\\

$\geq(3\delta-3\kappa+f-1)+2f-2m-p_0$.\\

Since $2f-2m-p_0\geq f-p_0\geq0$, we have $|V^{\downarrow}|\geq3\delta-3\kappa+f-1$.\\

\textbf{Case 1.2.1.2.2}. Either $|Z_{11}|\geq4$ or $|Z_{31}|\geq4$.

In this case we can argue exactly as in proof of Lemma 14 (Cases 2.2.1.2.2-2.2.1.2.3).\\

\textbf{Case 1.2.2}. $p_0\geq f+1$.

As in proof of Lemma 14 (Case 2.2.2), there are $f+1$ vertex disjoint paths $R_1,R_2,...,R_{f+1}$ connecting $V^{\downarrow}$ and $V(P_0)$. Let $F(R_i)\in V(P_0)$ and $L(R_i)\in V^{\downarrow}$ $(i=1,...,f+1)$. Since $f\geq3$ and $f+1\geq2m+1$, we can assume w.l.o.g. that either $L(R_i)\in V^{\downarrow}_1$ $(i=1,2,3,4)$ or $L(R_i)\in V^{\downarrow}_1$ $(i=1,2,3)$ and $L(R_i)\in V^{\downarrow}_2$ $(i=4,5)$.\\

\textbf{Case 1.2.2.1}. $L(R_i)\in V^{\downarrow}_1$ $(i=1,2,3,4)$.

Assume w.l.o.g. that $R_1\cup R_2$ dose not intersect $y_1\overrightarrow{P}w^-_1$ (see the proof of Lemma 12, Case 1.2.2.1).\\

\textbf{Case 1.2.2.1.1}. $\delta-\kappa+f-|Z_3|\leq1$.

Clearly $|Z_1|\geq|Z_3|\geq\delta-\kappa+f-1$. Using (12) and summing, \\

$|V^{\downarrow}|=\sum^{m}_{i=1}|V^{\downarrow}_i|\geq2(|Z_1|+|Z_3|)-2m\geq4(\delta-\kappa+f-1)-2m$\\

$=(3\delta-3\kappa+f-1)+\delta-\kappa+3f-2m-3>3\delta-3\kappa+f-1$.\\

\textbf{Case 1.2.2.1.2}. $\delta-\kappa+f-|Z_3|\geq2$.

By the choice of $w$, $|N(w_i)\cap V(P_0)|\geq\delta-\kappa+f-|Z_3|$ $(i=1,...,s)$. In particular, when $i=s$, we have $s\geq\delta-\kappa+f-|Z_3|$. By Lemma C, in $\langle V(P_0)\rangle$ any two vertices are joined by a path of length at least $\delta-\kappa+f-|Z_3|$. If $|Z_3|\leq3$, then due to $R_1,R_2,R_3,R_4$, \\

$|V^{\downarrow}|\geq|V^{\downarrow}_1|\geq3(\delta-\kappa+f-|Z_3|+2)+1$\\

$=(3\delta-3\kappa+f-1)+2f+8-3|Z_3|$.\\

If $|Z_3|\leq3$, then clearly we are done. Otherwise, we have $|Z_1|\geq|Z_3|\geq4$.\\

\textbf{Case 1.2.2.1.2.1}. $|Z_{11}|\leq3$, $|Z_{31}|\leq3$.

As in previous case, $|V^{\downarrow}_1|\geq3(\delta-\kappa+f-|Z_3|+2)+1$. Using (12) for each $i\in\{2,...,m\}$, we get\\

$|V^{\downarrow}|\geq3(\delta-\kappa+f-|Z_3|+2)+1+2(|Z_1|-|Z_{11}|+|Z_3|-|Z_{31}|)-2(m-1)$\\

$=(3\delta-3\kappa+f-1)+2|Z_1|-|Z_3|+2f-2m+10-2|Z_{11}|-2|Z_{31}|>3\delta-3\kappa+f-1$.\\

\textbf{Case 1.2.2.1.2.2}. $|Z_{11}|\geq4$, $|Z_{31}|\leq3$.

We can assume w.l.o.g. that $(Z_{11},\{L(R_1)\},\{L(R_2)\},\{L(R_3)\})$ is a nontrivial $(Q^{\downarrow}_1,\delta-\kappa+f-|Z_3|+2)$-scheme. By Lemma 8,\\

$|V^{\downarrow}_1|\geq3(\delta-\kappa+f-|Z_3|+2)+2|Z_{11}|-7$.\\

Using (11) for each $i\in \{2,...,m\}$ and summing, we get\\

$|V^{\downarrow}|\geq3(\delta-\kappa+f-|Z_3|+2)+2|Z_{11}|-7+2(|Z_1|-|Z_{11}|+|Z_3|-|Z_{31}|)-2(m-1)$\\

$=(3\delta-3\kappa+f-1)+(|Z_1|-|Z_3|)+|Z_1|+2f-3m+2-2|Z_{31}|>3\delta-3\kappa+f-1$.\\

\textbf{Case 1.2.2.1.2.3}. $|Z_{11}|\geq4$, $|Z_{31}|\geq4$.

We can assume w.l.o.g. that $(\{L(R_1)\},\{L(R_2)\},Z_{11},Z_{31})$ is a nontrivial $(Q^{\downarrow}_1,\delta-\kappa+f-|Z_3|+2)$-scheme. By Lemma 7,\\

$|V^{\downarrow}_1|\geq3(\delta-\kappa+f-|Z_3|)+2(|Z_{11}|+|Z_{31}|)-9$\\

$+\min\{2,\frac{1}{2}(\delta-\kappa+f-|Z_3|-2)(|Z_{11}|+|Z_{31}|-8)\}$\\

$\geq3(\delta-\kappa+f-|Z_3|)+2(|Z_{11}|+|Z_{31}|)-9$.\\

Using (12) for each $i\in\{2,...,m\}$ and summing, we get\\

$|V^{\downarrow}|\geq3(\delta-\kappa+f-|Z_3|)+2(|Z_{11}|+|Z_{31}|)-9$\\

$+2(|Z_1|-|Z_{11}|+|Z_3|-|Z_{31}|)-2(m-1)$\\

$=(3\delta-3\kappa+f-1)+2|Z_1|-|Z_3|+2f-2m-6\geq3\delta-3\kappa+f-1$.\\

\textbf{Case 1.2.2.2}. $L(R_i)\in V^{\downarrow}_1$ $(i=1,2,3)$ and $L(R_i)\in V^{\downarrow}_2$ $(i=4,5)$

In this case we can argue as in Case 1.2.2.1.\\ 

\textbf{Case 2}. $f=2$ and $S\not\subseteq V^{\uparrow}$.

Suppose first that $\delta-\kappa\leq1$. Since $\delta>3\kappa /2-1$, we have $|A^{\downarrow}|\geq\delta-\kappa+1\geq2$. Then it is easy to show that $|V^{\downarrow}|\geq4\geq3\delta-3\kappa+1$ and we are done. So, we can assume that $\delta-\kappa\geq2$. Let $S=\{v_1,...,v_{\kappa}\}$ and $F(Q^{\downarrow}_1)=v_1,L(Q^{\downarrow}_1)=v_2$. Assume w.l.o.g. that $v_3\in S-V^{\uparrow}$. Clearly $v_3\notin V^{\downarrow}$. Consider the graph $G^{\prime}=G-\{v_4,v_5,...,v_{\kappa}\}$. If there are two paths in $\langle A^{\downarrow}\cup \{v_1,v_2,v_3\}\rangle$ joining $v_3$ to $Q^{\downarrow}_1$ and having only $v_3$ in common, then the existence of $Q^{\downarrow}_{0}$ follows easily. Otherwise there is a cut-set in $\langle A^{\downarrow}\cup \{v_1,v_2,v_3\}\rangle$ consisting of a single vertex $z$ that separates $v_3$ and $V(Q^{\downarrow}_1)-z$. Then $\{v_1,v_2,z\}$ is an another cut-set of $G^{\prime}$, contradicting the definition of $A^{\downarrow}$. So, the existence of $Q^{\downarrow}_{0}$ is proved. By the definition of $Q^{\downarrow}_1$, we have $|Q^{\downarrow}_0|\leq |Q^{\downarrow}_1|$. The notation $P=y_1...y_p$, $Z_1$, $Z_2$ defined for $Q^{\downarrow}_1$, we will use here for $Q^{\downarrow}_0$. Let $M_1=v_1\overrightarrow{Q}^{\downarrow}_{0}v_3$ and $M_2=v_3\overrightarrow{Q}^{\downarrow}_{0}v_2$. In addition, put $Z_{1,i}=Z_1\cap V(M_i)$ and  $Z_{2,i}=Z_2\cap V(M_i)$ $(i=1,2)$. If $p\leq1$, then$\langle A^{\downarrow}-V(Q^{\downarrow}_0)\rangle$ is edgeless. Let $p\geq2$. If $v_3\in Z_1\cup Z_2$, then we can argue as in Case 1. Let $v_3\notin Z_1\cup Z_2$. Denote by $M^{\prime}_1$ and $M^{\prime}_2$ the minimal subsegments in $M_1$ and $M_2$, respectively, such that $Z_{1,1}\cup Z_{2,1}\subseteq V(M^{\prime}_1)$ and $Z_{1,2}\cup Z_{2,2}\subseteq V(M^{\prime}_2)$.\\ 

\textbf{Case 2.1.} $p=2$.

Clearly $|Z_i|\geq\delta-\kappa+2$ $(i=1,2)$. Applying (6) to $M_1^{\prime}$ and $M_2^{\prime}$, we get\\

$|V^{\downarrow}|=|V^{\downarrow}_{1}|\geq |V(Q^{\downarrow}_0)|\geq|V(M_1^{\prime})|+|V(M_2^{\prime})|+1$\\

$\geq 3(|Z_1|+|Z_2|)/2-3\geq3(\delta-\kappa+2)-3>3\delta-3\kappa+1.$\\

\textbf{Case 2.2.} $p=3$.

Clearly $|Z_i|\geq\delta-\kappa+1$ $(i=1,2)$.  Applying (8) to $M^{\prime}_1$ and $M^{\prime}_2$ , we get\\ 

$|V^{\downarrow}|\geq |V(Q^{\downarrow}_0)|\geq|V(M_1^{\prime})|+|V(M_2^{\prime})|+1$\\

$\geq 2(|Z_1|+|Z_2|)-5\geq4(\delta-\kappa+1)-5\geq3\delta-3\kappa+1.$\\

\textbf{Case 2.3.} $p\geq4$.

Let $P_0,p_0,w,Z_3$ are as defined in Lemma 14 (Case 2). If $p_0\leq3$, then $|Z_1|\geq|Z_2|\geq\delta-\kappa+1$ and we can argue as in Case 2.2. Let $p_0\geq4$. Further, if $\delta-\kappa+3-|Z_3|\leq1$, then $|Z_1|\geq|Z_3|\geq\delta-\kappa+2$ and we can argue as in Case 2.1. Let $\delta-\kappa+3-|Z_3|\geq2$. Since $p_0\geq4$, there are vertex disjoint paths $R_1,R_2,R_3,R_4$ in $G^{\prime}$ connecting $P_0$ and $Q^{\downarrow}_{0}$ (otherwise, there exist a cut-set of $G^{\prime}$ of order 3 contradicting the definition of $A^{\downarrow}$). Let $F(R_i)\in V(P)$ and $L(R_i)\in V(Q^{\downarrow}_{0})$ $(i=1,2,3,4)$. Clearly, $|N(w_i)\cap V(P_0)|\geq\delta-\kappa+3-|Z_3|$ for each $i\in\{1,...,s\}$. In particular, for $i=s$, we have $s\geq\delta-\kappa+3-|Z_3|$. By Lemma C, in $\langle V(P_0)\rangle$ any two vertices are joined by a path of length at least $\delta-\kappa+3-|Z_3|$. Assume w.l.o.g. that $R_1\cup R_2$ does not intersect $y_1\overrightarrow{P}w^{-}_{1}$ (see the proof of Lemma 12, Case 1.2.2.1) and does not contain $w$. Let $I_1,...,I_t$ be the minimal segments of $Q^{\downarrow}_{0}$ connecting two vertices of $Z_1\cup Z_3\cup Z_4\cup Z_5$, where $Z_4=\{L(R_1)\}$ and $Z_5=\{L(R_2)\}$. Assume w.l.o.g. that $v_3$ belongs to $I_1$. Put $I_1=v^{\prime}\overrightarrow{Q}^{\downarrow}_{1}v^{\prime\prime}$. If $v_3=v^{\prime}$ or $v_3=v^{\prime\prime}$, then we can argue as in Case 1. Let $v_3\neq v^{\prime}$ and $v_3\neq v^{\prime\prime}$. Choose a longest path $Q_0$ joining $v^{\prime}$ and $v^{\prime\prime}$ and passing through $V(P_0)\cup^{4}_{i=1}V(R_i)$. Clearly, $|Q_0|\geq\delta-\kappa+5-|Z_3|$ if $v^{\prime},v^{\prime\prime}$ belong to different $Z_1,Z_3,Z_4,Z_5$ and $|Q_0|\geq2$, otherwise. Since $Q^{\downarrow}_{0}$ is extreme with ends $v_1,v_2$ and intermediate vertex $v_3$, we have $|I_i|\geq\delta-\kappa+5-|Z_3|$ if the ends of $I_i$ belong to different $Z_1,Z_3,Z_4,Z_5$ and $|I_i|\geq2$ for each $i\in\{2,...,t\}$, otherwise. Form a new path $Q_{00}$ from $Q^{\downarrow}_{0}$ by replacing $I_1$ with $Q_0$. Clearly $(Z_1,Z_3,Z_4,Z_5)$ is a nontrivial $(Q_{00},\delta-\kappa+5-|Z_3|)$-scheme and we can argue as in Case 1.2.2.1. \\

\textbf{Case 3.} $f=2$ and $S\subseteq V^{\uparrow}$.

Let$P_0$, $p_0$, $w$ and $Z_3$ are as defined in Lemma 14 (Case 2). If $p_0\leq2$, then $|Z_1|\geq|Z_2|\geq\delta-\kappa+1$ and we can argue as in proof of Lemma 14 (Case 1). Let $p_0\geq3$. Since $\kappa\geq3$, there are vertex disjoint paths $R_1,R_2,R_3$ connecting $V^{\downarrow}$ and $V(P_0)$ (see the proof of Lemma 14, Case 2.2.2). Let $F(R_i)\in V(P_0)$ and $L(R_i)\in V^{\downarrow}$ $(i=1,2,3)$. If $\delta-\kappa-|Z_3|+2\leq1$, then $|Z_1|\geq|Z_3|\geq\delta-\kappa+1$ and again we can argue as in the case $p_0\leq2$. Let $\delta-\kappa-|Z_3|+2\geq2$. By the choice of $w$, $|N(w^-_i)\cap V(P_0)|\geq\delta-\kappa+f-|Z_3|$ $(i=1,...,s)$. In particular, for $i=s$, $s\geq\delta-\kappa-|Z_3|+2$. By Lemma C, in $\langle V(P_0)\rangle$ any two vertices are joined by a path of length at least $\delta-\kappa-|Z_3|+2$. If $|Z_1|\leq2$, then $|Z_3|\leq2$. Since $(\{ L(R_1)\},\{L(R_2)\},\{L(R_3)\})$ is a nontrivial $(Q^{\downarrow}_1,\delta-\kappa-|Z_3|+4)$-scheme, we have $|V^{\downarrow}|\geq2(\delta-\kappa-|Z_3|+4)+1\geq2(\delta-\kappa+2)+1\geq2\delta-2\kappa+f+1$. Let $|Z_1|\geq3$. Analogously, $|Z_3|\geq3$. Assume w.l.o.g. that $F(R_1)=w$. In addition, we can assume w.l.o.g. that $R_2\cup R_3$ does not intersect $y_1\overrightarrow{P}w^-_1$ (see the proof of Lemma 12, Case 1.2.2.1). So, $(Z_1,Z_3,\{L(R_1)\})$ is a nontrivial $(Q^{\downarrow}_1,\delta-\kappa-|Z_3|+4)$-scheme and by Lemma 5,\\

$|V^{\downarrow}|\geq\min\{2(|Z_1|+|Z_3|+1)+2(\delta-\kappa-|Z_3|+4)-11,$\\

$\frac{1}{2}(\delta-\kappa-|Z_3|+4)(|Z_1|+|Z_3|-2)+1\}$\\

$\geq\min\{2(2|Z_3|+1)+2(\delta-\kappa)-2|Z_3|-3, (\delta-\kappa-|Z_3|+4)(|Z_3|-1)+1\}$\\

$\geq2\delta-2\kappa+2|Z_3|-3+\min\{2,(\delta-\kappa-|Z_3|)(|Z_3|-3)\}\geq2\delta-2\kappa+f+1$.\\

Lemma 15 is proved. \quad $\Delta$

\section{Proofs of theorems}

\textbf{Proof of Theorem 2}. If $\delta\leq3k/2-1$, then we are done by Lemma 11. Let $\delta>3\kappa/2-1$. We will use the notation defined in Definitions A and B. In addition, let $A^{\uparrow}$ is defined with respect to a minimum cut-set $S=\{v_1,...,v_{\kappa}\}$. The existence of $Q^{\downarrow}_1,...,Q^{\downarrow}_m$ and $C^{*},C^{**}$ follows from Lemma 10. By Lemmas 12 and 14, $A^{\uparrow}\subseteq V^{\uparrow}$ and $\langle A^{\downarrow}-V^{\downarrow}\rangle$ is edgeless. Since $Q^{\uparrow}_1,...,Q^{\uparrow}_m$ is extreme, $\langle S-V(C^{*})\rangle$ is edgeless too. Recalling the definition of $C^{**}$, we can state that $A^{\uparrow}\subseteq V(C^{**})$ and in addition, $A^{\downarrow}-V(C^{**})$ and $S-V(C^{**})$ both are edgeless.\\

\textbf{Case 1}. $A^{\downarrow}-V(C^{**})=\emptyset$.

If $S-V(C^{**})=\emptyset$, then $C^{**}$ is a Hamilton Cycle. Let $S-V(C^{**})\neq\emptyset$ and choose any $w\in S-V(C^{**})$. Clearly $N(w)\subseteq V(C^{**})$. If $N(w)^+\cup \{w\}$ is independent, then $\alpha\geq\delta+1$, contradicting the hypothesis. Otherwise we can form (by standard arguments) a cycle with vertex set $V(C^{**})\cup\{w\}$, contradicting the definition of $C^{**}$.\\

\textbf{Case 2}. $A^{\downarrow}-V(C^{**})\neq\emptyset$.

Let $z\in A^{\downarrow}-V(C^{**})$. If $N(z)\subseteq V(C^{**})$, then we can argue as in Case 1. Otherwise $N(z)=D_1\cup D_2$, where $D_1\subseteq V(C^{**})$ and $D_2\subseteq S-V(C^{**})$. If $|A^{\uparrow}|\leq\kappa$, then $n\leq|A^{\uparrow}|+|A^{\downarrow}|+\kappa\leq3\delta-\kappa$ and by Theorem B, $G$ is hamiltonian. Let $|A^{\uparrow}|\geq\kappa+1$. If $N(v_i)\cap A^{\uparrow}=\emptyset$ for some $i\in\{1,...,\kappa\}$, then $S-v_i$ is a cut-set of order $\kappa-1$, a contradiction. Therefore, 

$$
N(v_i)\cap A^{\uparrow}\neq\emptyset \quad (i=1,...,\kappa).\eqno{(13)} 
$$

Put $N_i=N(v_i)\cap A^{\uparrow} \quad (i=1,...,\kappa)$. If $|\cup_{i\in J}N_i|<|J|\leq\kappa$ for a subset $J\subseteq\{1,...,\kappa\}$, then $(\cup_{i\in J}N_i)\cup\{v_i\in S|i\notin J\}$ is a cut-set of $G$ with at most $\kappa-1$ vertices, a contradiction. So, we can assume that $|\cup_{i\in J}N_i|\geq |J|$ for each $J\subseteq\{1,...,\kappa\}$. By Hall's [5] Theorem, the collection $N_1,...,N_{\kappa}$ has a system of distinct representatives. Set $D_2=\{v_{{i}_1},v_{{i}_2},...,v_{{i}_t}\}$ and let $w_{{i}_1},...,w_{{i}_t}$ is a system of distinct representatives of $N_{{i}_1},...,N_{{i}_t}$ . Put $D_3=\{w_{{i}_1},...,w_{{i}_t}\}$. Since $A^{\uparrow}\subseteq V(C^{**})$, we have $D_3\subseteq V(C^{**})$ and it is easy to see that $(D_1\cup D_3)^+\cup\{z\}$ is an independent set of order at least $\delta+1$, a contradiction. \quad $\Delta$\\

\noindent\textbf{Proof of Theorem 3}. If $\delta\leq3k/2-1$, then we are done by Lemma 11. Let $\delta>3\kappa/2-1$. By Lemma 12, $A^{\uparrow}\subseteq V^{\uparrow}$. The existence of $Q^{\downarrow}_1,...,Q^{\downarrow}_m$ and $C^*,C^{**}$ (see Definition B) follows from Lemma 10. Let $A^{\downarrow}$ is defined with respect to a minimum cut-set $S=\{v_1,...,v_{\kappa}\}$. Put $f=|V^{\downarrow}\cap S|$. By Theorem G, $c\geq\min\{n,3\delta-3\}$. If $c\geq n$, then we are done. Let $c\geq 3\delta-3$. Further, if $3\delta-3\geq4\delta-2\kappa$, i.e. $\delta\leq2\kappa-3$, then $c\geq3\delta-3\geq 4\delta-2\kappa$. Let $\delta\geq2\kappa-2$ which implies $|A^{\downarrow}|\geq3\delta-3\kappa+1\geq2\delta-\kappa-1$. Recalling also that $|A^{\uparrow}|\geq|A^{\downarrow}|$, we obtain

$$
\delta-2\kappa+2\geq0, \quad |A^{\uparrow}|\geq|A^{\downarrow}|\geq2\delta-\kappa-1. \eqno{(14)}
$$

If $A^{\downarrow}\subseteq V^{\downarrow}$, then by (14), $c\geq |A^{\uparrow}|+|A^{\downarrow}|+2\geq 4\delta-2\kappa$. Let

$$
A^{\downarrow}\not\subseteq V^{\downarrow}. \eqno{(15)}
$$

\textbf{Case 1}. $f\geq3$.

By Lemma 15, either $\langle A^{\downarrow}-V^{\downarrow}\rangle$ is edgeless or $|V^{\downarrow}|\geq3\delta-3\kappa+f-1$. If $\langle A^{\downarrow}-V^{\downarrow}\rangle$ is edgeless, then we can argue as in proof of Theorem 2. Let $|V^{\downarrow}|\geq3\delta-3\kappa+f-1$. By (14), \\

$c\geq|A^{\uparrow}|+|V^{\downarrow}|\geq(2\delta-\kappa-1)+(3\delta-3\kappa+f-1)$\\

$=(4\delta-2\kappa)+(\delta-2\kappa+2)+f-4\geq4\delta-2\kappa+f-4$.\\

If $f\geq4$, then we are done. Let $f=3$. Then $c\geq|A^{\uparrow}|+|V^{\uparrow}\cap S|-2+|V^{\downarrow}|$ and the desired result can be obtained by a similar calculation as above, if either $|A^{\uparrow}|\geq2\delta-\kappa$ or $|V^{\uparrow}\cap S|\geq3$. So, we can assume that $|A^{\uparrow}|=|A^{\downarrow}|=2\delta-\kappa-1$ and $|V^{\uparrow}\cap S|=2$. Let $V^{\uparrow}\cap S=\{v_1,v_2\}$ and $V^{\downarrow}\cap S=\{v_1,v_2,v_3\}$. Consider the graph $G^{\prime}=G-\{v_4,...,v_{\kappa}\}$. By (15), there exists a connected component $H$ of $G^{\prime}-Q^{\downarrow}_{1}$ intersecting $A^{\downarrow}$. Put $M=N(V(H))$ and let $x_1,...,x_q$ be the elements of $M\cap V^{\downarrow}$ occuring on $\overrightarrow{Q}^{\downarrow}_{1}$ in a consecutive order. Since $G^{\prime}$ is 3-connected, we have $q\geq 3$. Further, $|x_i\overrightarrow{Q}^{\downarrow}_{1}x_{i+1}|\geq2$ $(i=1,...,q-1)$, since $Q^{\downarrow}_{1}$ is extreme. Put\\

$M^*=\{x\in M\cap V^{\downarrow}|x^-\notin S\quad \mbox{or}\quad x^+\notin S\}$.\\

If $M^{*}=\emptyset$, then it is easy to check that $|V^{\downarrow}\cap S|\geq4$, contradicting the fact that $f=3$. Let $M^{*}\neq\emptyset$ and choose any $u\in M^{*}$. Assume w.l.o.g. that $u^+\notin S$. Choose $w\in V(H)$ such that $uw\in E(G)$. Then by deleting $uu^+$ from $Q^{\downarrow}_{1}$ and adding $uw$, we obtain a pair of vertex disjoint paths $v_1\overrightarrow{Q}^{\downarrow}_{1}w$ and $v_2\overleftarrow{Q}^{\downarrow}_{1}u^+$ both having one end in $S$ and the other in $A^{\downarrow}$. These two paths can be extended from $w$ and $u^+$ (using only vertices of $A^{\downarrow}$) to a pair of maximal vertex disjoint paths $R_1=v_1\overrightarrow{R}_1w_1$ and $R_2=v_2\overrightarrow{R}_2w_2$. Add an extra edge $v_1v_2$ to $G$ and consider the path $L=w_1\overleftarrow{R}_1v_1v_2\overrightarrow{R}_2w_2$ in $G^*$. Label $L=\xi_1\xi_2...\xi_h$ according to the direction on $L$ and put \\

$d_1=|N(\xi_1)\cap V(L)|,\quad d_2=|N(\xi_h)\cap V(L)|$. \\

If $\xi_1$ and $\xi_h$ have a common neighbor $v_i$ in $\{v_4,...,v_{\kappa}\}$, then $v_1\overrightarrow{R}_1\xi_1v_i\xi_h\overleftarrow{R}_2v_2$ is a path contradicting the choice of $Q^{\downarrow}_{1}$. Otherwise, $d_1+d_2\geq 2\delta-\kappa+3$. Since $|V(L)|\leq |A^{\downarrow}|+3=2\delta-\kappa+2$, i.e. $d_1+d_2\geq |V(L)|+1$, it can be shown by standard arguments that $\xi_1\xi_{i+1}\in E(G)$ and $\xi_h\xi_i\in E(G)$ for some distinct $i=i_1,i_2$. Assume w.l.o.g. that $\xi_{i_1}\xi_{i_1+1}\neq v_1v_2$ and $\xi_{i_1}\in \xi_1\overrightarrow{L}v_1$. Then $v_1\overleftarrow{L}\xi_{i_1+1}\xi_1\overrightarrow{L}\xi_{i_1}\xi_h\overleftarrow{L}v_2$ is a path longer than $Q^{\downarrow}_{1}$ contradicting the definition of $Q^{\downarrow}_{1}$. \\

\textbf{Case 2}. $f=2$ and $S\subseteq V^{\uparrow}$.

By Lemma 15, either $\langle A^{\downarrow}-V^{\downarrow}\rangle$ is edgeless or $|V^{\downarrow}|\geq2\delta-2\kappa+3$. In the first case we can argue as in proof of Theorem 2. In the second case,\\

$c\geq|A^{\uparrow}|+|S|+|V^{\downarrow}|-f\geq(2\delta-\kappa-1)+\kappa+(2\delta-2\kappa+3)-2=4\delta-2\kappa$.\\

\textbf{Case 3}. $f=2$ and $S\not\subseteq V^{\uparrow}$.

By Lemma 15, either $\langle A^{\downarrow}-V(Q^{\downarrow}_0)\rangle$ is edgeless or $|V^{\downarrow}|\geq3\delta-3\kappa+1$.\\

\textbf{Case 3.1}. $\langle A^{\downarrow}-V(Q^{\downarrow}_0)\rangle$ is edgeless.

By the definition of $C^{**}_0$, we have $A^{\uparrow}\subseteq V(C^{**}_0)$. Besides, $A^{\downarrow}-V(C^{**}_0)$ and $S-V(C^{**}_0)$ both are independent in $G$. If $A^{\downarrow}-V(C^{**}_0)\neq\emptyset$, then we can argue as in proof of Theorem 2 (Case 2). Otherwise \\

$c\geq|A^{\uparrow}|+|A^{\downarrow}|+3\geq2(2\delta-\kappa-1)+3>4\delta-2\kappa$.\\ 

\textbf{Case 3.2}. $|V^{\downarrow}|\geq3\delta-3\kappa+1$.

Let $R_1=v_1\overrightarrow{R}_1w_1$, $R_2=v_2\overrightarrow{R}_2w_2$, $L=\xi_1...\xi_h$ and $d_1,d_2$ be as defined in Case 1 with respect to $Q^{\downarrow}_1$. Put $|V^{\uparrow}\cap S|=f^{\prime}$. Using (14), we get\\

$$c\geq |A^{\uparrow}|+|V^{\downarrow}|+f^{\prime}-2\geq2\delta-\kappa-1+3\delta-3\kappa+1+f^{\prime}-2$$

$$=(4\delta-2\kappa)+\delta-2\kappa+2+f^{\prime}-4\geq4\delta-2\kappa+f^{\prime}-4. \eqno{(16)}$$

If $f^{\prime}\geq4$, then we are done. Let $f^{\prime}\leq3$. Similar to (16), we can state that

$$
\mbox{if} \quad|A^{\uparrow}|\geq2\delta-\kappa, \quad \mbox{then} \quad c\geq4\delta-2\kappa+f^{\prime}-3, \eqno{(17)}\\
$$

$$
\mbox{if} \quad|A^{\uparrow}|\geq2\delta-\kappa+1, \quad \mbox{then} \quad c\geq4\delta-2\kappa \eqno{(18)}\\
$$

\textbf{Case 3.2.1}. $f^{\prime}=3$. 

If $|A^{\uparrow}|\geq2\delta-\kappa$, then by (17) we are done. Let $|A^{\uparrow}|=2\delta-\kappa-1$, implying also $|A^{\downarrow}|=2\delta-\kappa-1$. If $\xi_1$ and $\xi_h$ have a common neighbor $v_i$ in $\{v_4,...,v_{\kappa}\}$, then $v_1\overrightarrow{R}_1\xi_1v_i\xi_h\overleftarrow{R}_2v_2$ is a path contradicting the choice of $Q^{\downarrow}_1$. Otherwise we have $d_1+d_2\geq2\delta-\kappa+1$. In addition, $|V(L)|\leq |A^{\downarrow}|+2=2\delta-\kappa+1$. If $|V(L)|<2\delta-\kappa+1$, then as in Case 1, $d_1+d_2\geq|V(L)|+1$ and we can form a path longer than $Q^{\downarrow}_1$, connecting $v_1,v_2$ and passing through $A^{\downarrow}$, contrary to the definition of $Q^{\downarrow}_1$. Hence, $|V(L)|=2\delta-\kappa+1$. This means that $v_3$ is adjacent to both $\xi_1$ and $\xi_h$. Besides, each $v_i$ $(i\in\{4,...,\kappa\})$ is adjacent either to $\xi_1$ or $\xi_h$. Assume w.l.o.g. that $v_4\xi_h\in E(G)$. Put $Y_1=v_1\overleftarrow{L}\xi_1v_3$ and $Y_2=v_2\overrightarrow{L}\xi_hv_4$. Since $|A^{\uparrow}|=|A^{\downarrow}|$, we can state that $A^{\uparrow}$ and $A^{\downarrow}$ are both endfragments. Then taking $\{Y_1,Y_2\}$ instead of $\{Q^{\uparrow}_1,...,Q^{\uparrow}_m\}$ and $A^{\downarrow}$ instead of $A^{\uparrow}$, we can argue as in case $f^{\prime}\geq4$.\\

\textbf{Case 3.2.2}. $f^{\prime}=2$.

If $|A^{\uparrow}|\geq2\delta-\kappa+1$, then we are done by (18). Let $|A^{\uparrow}|\leq2\delta-\kappa$ implying also $|A^{\downarrow}|\leq2\delta-\kappa$. Further, we have $d_1+d_2\geq2\delta-\kappa+2$.\\

\textbf{Case 3.2.2.1}. $|A^{\uparrow}|=2\delta-\kappa-1$.

In this case, $|A^{\downarrow}|=2\delta-\kappa-1$. Clearly $|V(L)|\leq|A^{\downarrow}|+2=2\delta-\kappa+1$ implying that $d_1+d_2\geq|V(L)|+1$. Then we can form (as above) a path contradicting the definition of $Q^{\downarrow}_1$.\\

\textbf{Case 3.2.2.2}. $|A^{\uparrow}|=2\delta-\kappa$.

In this case, $2\delta-\kappa-1\leq|A^{\downarrow}|\leq2\delta-\kappa$. If $|A^{\downarrow}|=2\delta-\kappa-1$, then $|V(L)|\leq|A^{\downarrow}|+2=2\delta-\kappa+1$ and hence $d_1+d_2\geq|V(L)|+1$. Then again we can form a path contradicting the choice of $Q^{\downarrow}_1$. Now let $|A^{\downarrow}|=2\delta-\kappa$. It follows that $A^{\uparrow}$ and $A^{\downarrow}$ are both endfragments. On the other hand, $|V(L)|\leq|A^{\downarrow}|+2=2\delta-\kappa+2$ implying that $d_1+d_2\geq|V(L)|$. Thus we can argue as in Case 3.2.1. \quad $\Delta$\\

\noindent\textbf{Proof of Theorem 4}. If $\delta\leq3k/2-1$, then we are done by Lemma 11. Let $\delta>3\kappa/2-1$. By Lemma 14, $\langle A^{\downarrow}-V^{\downarrow}\rangle$ is edgeless. \\

\textbf{Case 1}. $A^{\downarrow}\subseteq V^{\downarrow}$.

If $A^{\uparrow}\subseteq V^{\uparrow}$, then \\

$c\geq|A^{\uparrow}|+|A^{\downarrow}|+2\geq(3\delta-\kappa-3)+(\delta-\kappa+1)+2=4\delta-2\kappa$.\\ 

Let $A^{\uparrow}\not\subseteq V^{\uparrow}$. By Lemma 13, $|V^{\uparrow}|\geq3\delta-5$ and hence\\

$c\geq|V^{\uparrow}|+|A^{\downarrow}|\geq3\delta-5+\delta-\kappa+1\geq 4\delta-2\kappa$.\\

\textbf{Case 2}. $A^{\downarrow}\not\subseteq V^{\downarrow}$.
 
By the definition of $C^{**}$, $\langle A^{\downarrow}-V(C^{**})\rangle$ is edgeless and hence $N(z)\subseteq V(C^{**})\cup S$ for each $z\in A^{\downarrow}-V(C^{**})$. If $N(z)\subseteq V(C^{**})$, then by standard arguments, $\alpha\geq\delta+1$, a contradiction. Let $N(z)=D_1\cup D_2$, where $D_1\subseteq V(C^{**})$ and $D_2\subseteq S-V(C^{**})$. Set $D_2=\{v_{i_1},...,v_{i_t}\}$ and $N_i=N(v_i)\cap A^{\uparrow}$ $(i=i_1,...,i_t)$. As in proof of Theorem 2, the collection $N_{i_1},...,N_{i_t}$ has a system of distinct representatives $w_{i_1},...,w_{i_t}$. Put $D_3=\{w_{i_1},...,w_{i_t}\}$. Let $D_3=D_4\cup D_5$, where $D_4\subseteq V(C^{**})$ and $D_5=D_3-D_4$. By the definition of $Q^{\uparrow}_{1},...,Q^{\uparrow}_{m}$, it is easy to see that $(D_1\cup D_4)^+\cup D_3\cup \{z\}$ is an independent set with at least $\delta+1$ vertices, contradicting $\delta\geq\alpha$. \quad $\Delta$\\

\noindent\textbf{Proof of Theorem 5}. If $\delta\leq3k/2-1$, then we are done by Lemma 11. Let $\delta>3\kappa/2-1$. The existence of $Q^{\downarrow}_1,...,Q^{\downarrow}_m$ and $C^{*},C^{**}$ follows from Lemma 10. As in proof of theorem 3, $\delta-2\kappa+2\geq0$, implying in particular that $\delta-\kappa\geq2$. By Lemma 13, either $A^{\uparrow}\subseteq V^{\uparrow}$ or $|V^{\uparrow}|\geq3\delta-5$. \\

\textbf{Case 1}. $A^{\downarrow}\subseteq V^{\downarrow}$. 

If $A^{\uparrow}\subseteq V^{\uparrow}$, then\\ 

$c\geq|A^{\uparrow}|+|V^{\downarrow}|+2\geq(3\delta-\kappa-3)+(3\delta-3\kappa+1)+2>4\delta-2\kappa$.\\

If $|V^{\uparrow}|\geq3\delta-5$, then \\

$c\geq|V^{\uparrow}|+|A^{\downarrow}|\geq3\delta-5+3\delta-3\kappa+1\geq4\delta-2\kappa$. \\

\textbf{Case 2}. $A^{\downarrow}\not\subseteq V^{\downarrow}$.

If either $\langle A^{\downarrow}-V^{\downarrow}\rangle$ or $\langle A^{\downarrow}-V(Q^{\downarrow}_0)\rangle$ is edgeless, then we can argue as in proof of Theorem 4. Otherwise, by Lemma 15, either $f=2$ and $|V^{\downarrow}|\geq2\delta-2\kappa+3$ or $f\geq3$ and $|V^{\downarrow}|\geq3\delta-3\kappa+f-1\geq2\delta-2\kappa+3$, where $f=|V^{\downarrow}\cap S|$. If $A^{\uparrow}\subseteq V^{\uparrow}$, then\\

$c\geq|A^{\uparrow}|+|V^{\downarrow}|\geq3\delta-\kappa-3+2\delta-2\kappa+3\geq4\delta-2\kappa$. \\

Let $A^{\uparrow}\not\subseteq V^{\uparrow}$. By Lemma 13, $|V^{\uparrow}|\geq3\delta-5$. If $f=2$, then\\

$c\geq|V^{\uparrow}|+|V^{\downarrow}|-2\geq(3\delta-5)+(2\delta-2\kappa+3)-2$\\

$=(4\delta-2\kappa)+\delta-4\geq 4\delta-2\kappa$. \\

If $f\geq3$, then\\

$c\geq|V^{\uparrow}|+|V^{\downarrow}|-f\geq(3\delta-5)+3\delta-3\kappa+f-1-f$.\\

$=(4\delta-2\kappa)+(\delta-2\kappa+2)+(\delta+\kappa-8)\geq4\delta-2\kappa$. \quad $\Delta$\\

\noindent\textbf{Proof of Theorem 1}. If $\delta\leq3\kappa/2-1$, then we are done by Lemma 11. Let $\delta>3\kappa/2-1$. Let $A^{\downarrow}$ be an endfragment of $G$ with $|A^{\uparrow}|\geq|A^{\downarrow}|$. Then the desired result follows from Theorems 2-5. \quad $\Delta$\\
 
\noindent\textbf{Remark}. The limit examples below show that Theorem 1 is best possible in all respects. The limit example $4K_2+K_3$ shows that 4-connectedness can not be replaced by 3-connectedness. Further, the limit example $H(1,1,5,4)$ shows that the condition $\delta\geq\alpha$ cannot be replaced by $\delta\geq\alpha-1$. Finally, the limit example $5K_2+K_4$ shows that the bound $4\delta-2\kappa$ cannot be replaced by $4\delta-2\kappa+1$.  

\section{Acknowledgments}

The author wishes to thank V.V. Mkrtchyan for his careful reading and his corrections for the manuscript.

\end{document}